\documentclass[preprint,3p,sort&compress,final,times]{elsarticle}
\usepackage{amsmath}
\usepackage{amssymb}
\usepackage{mathtools}
\usepackage{placeins}
\usepackage[usenames,dvipsnames]{color}
\usepackage{xcolor}
\usepackage{cases}
\usepackage{array}
\usepackage{todonotes}
\usepackage{lineno}
\newcommand{\blue}[1]{\color{black}{#1}\color{black}}

\journal{}
\begin{document}
\begin{frontmatter}

\title{Entropy and energy conservation for thermal atmospheric dynamics using mixed compatible finite elements}
\author[ANU]{Kieran Ricardo}
\author[BOM]{David Lee\corref{cor}}
\ead{davelee2804@gmail.com}
\author[ANU]{Kenneth Duru}

\address[ANU]{Mathematical Sciences Institute, Australian National University, Canberra, Australia}
\address[BOM]{Bureau of Meteorology, Melbourne, Australia}
\cortext[cor]{Corresponding author.}

\begin{abstract}
Atmospheric systems incorporating thermal dynamics must be stable with respect to both energy and entropy. 
While energy conservation can be enforced via the preservation of the skew-symmetric structure of the 
Hamiltonian form of the equations of motion, entropy conservation is typically derived as an additional 
invariant of the Hamiltonian system, and satisfied via the exact preservation of the chain rule. This is 
particularly challenging since the function spaces used to represent the thermodynamic variables in 
compatible finite element discretisations are typically discontinuous at element boundaries.
In the present work we negate this problem by constructing our equations of 
motion via a novel formulation
which allow for the necessary cancellations required to simultaneously conserve entropy and energy
without the chain rule. We show that such formulations allow for stable simulation of turbulent
dynamics for both the thermal 
shallow water and 3D compressible Euler equations on the sphere using mixed compatible finite elements 
without entropy damping.
\end{abstract}

\end{frontmatter}

Entropy stability is an essential feature of numerical schemes for the solution of hyperbolic systems 
involving thermally forced dynamics. In recent years great advances have been made in the formulation of 
entropy stable schemes using finite difference \cite{FC13}, spectral element \cite{Carpenter14} and 
discontinuous Galerkin \cite{GWK16,Chan18} methods for the solution of hyperbolic systems, including the 3D
compressible Euler equations used for the simulation of atmospheric dynamics \cite{Waru22}. These
formulations are based on the preservation of the summation by parts (SBP) property between gradient and
divergence operators for flux form representations of these systems. In a separate line of research, 
a finite volume scheme for the rotating thermal shallow water equations has been developed using global 
integrals of the source terms in order to derive well balanced fluxes that preserve equilibrium solutions 
\cite{KLZ20}.

One draw back of these flux 
form energy and entropy conserving schemes is that they do not naturally conserve vorticity, as do 
vector invariant formulations using compatible finite element schemes \cite{MC14,LPG18,Eldred19}.
More recently, an entropy stable discontinuous Galerkin formulation has been developed for the 
vector invariant thermal shallow water equations \cite{Ricardo23b} which builds off previous work on 
the shallow water equations \cite{Ricardo23a} to incorporate the properties of geostrophic balance 
and vorticity conservation, making this an appealing scheme for the simulation of atmospheric 
dynamics. In addition to the SBP property, this vector invariant formulation negates 
the need to preserve the chain rule in order to conserve entropy for the case of thermodynamic 
variables represented by discontinuous function spaces, for which the chain rule is not defined at 
element boundaries.

In the present work we extend this approach to the mixed finite element method, using compatible
function spaces. Moreover we show that this entropy conserving formulation
bears a strong similarity to an averaging
of two different skew-symmetric Hamiltonian formulations, one involving the material form transport
of the thermodynamic variable, and the other using flux form transport.
In doing so we show that thermal atmospheric systems including the thermal shallow water equations
and the 3D compressible Euler equations can run stably without any entropy dissipation, as is typically
required for the stable simulation of these systems. Note that in the present article we are referring
specifically to a \emph{mathematical} entropy, by which we mean a conserved moment that remains convex
for all possible solutions due the presence of strictly positive eigenvalues for its Hessian matrix 
\cite{FC13}, and not a \emph{physical} entropy, as may be derived from the functional derivatives of 
a thermodynamic potential \cite{Thuburn17}.

The remainder of this article proceeds as follows:
In Section 1 we will introduce the thermal shallow water and 3D compressible Euler equations in continuous
form, as well as their entropy functions. In Section 2 we describe discrete, mixed finite element 
formulations of these systems for which energy and entropy are conserved. We motivate these formulations
via an averaging of discrete skew-symmetric Hamiltonian systems, and then extend this approach to a coupled 
system using flux form transport only, which is similar to the approach previously derived for
the discontinuous Galerkin method \cite{Ricardo23b}. Results will be presented in Section 3, and finally
we will discuss the conclusions from this work in Section 4.

\section{Atmospheric systems with thermal dynamics}

\subsection{Thermal shallow water equations}\label{sec::tsw_cont}

The thermal shallow water equations may be considered as an extension of the shallow water equations
for which the buoyancy, $b$ is allowed to vary with respect to changes in density, $\rho$ as 
$b=g\rho/\overline{\rho}$, where $g$ is gravity and $\overline{\rho}$ is the mean density. If there is
no variation in density, such that $\rho=\overline{\rho}$, then the original shallow water equations
are recovered for a constant buoyancy. The Hamiltonian structure of the thermal shallow water equations
has been studied in continuous form in \cite{Dellar03,WD13}, and in a discrete mixed Galerkin variational 
form in \cite{Eldred19}. These may be given for the velocity, $\boldsymbol{u}$, fluid depth, $h$ and 
buoyancy $b$ as

\begin{subequations}
\begin{align}
	\frac{\partial\boldsymbol{u}}{\partial t} + q\times\boldsymbol{F} + 
	\nabla\Bigg(\frac{1}{2}\boldsymbol{u}\cdot\boldsymbol{u} + hb\Bigg) -
	\frac{h}{2}\nabla b &= 0\\
	\frac{\partial h}{\partial t} + \nabla\cdot\boldsymbol{F} &= 0\label{eq::tsw_cont}\\
	\frac{\partial b}{\partial t} + \frac{\nabla b}{h}\cdot\boldsymbol{F} &= 0,\label{eq::tsw_mat_adv}
\end{align}
\end{subequations}
where $\boldsymbol{F} = h\boldsymbol{u}$ is the mass flux and $q = (\nabla\times\boldsymbol{u} + f)/h$
is the potential vorticity (with $f$ being the Coriolis term due to the earth's rotation).
Combining \eqref{eq::tsw_cont} and \eqref{eq::tsw_mat_adv} we may alternatively express the buoyancy
transport as a flux form equation for the density weighted buoyancy $B=hb$ as
\begin{equation}
	\frac{\partial B}{\partial t} + \nabla\cdot(b\boldsymbol{F}) = 0.\label{eq::tsw_flux_adv}
\end{equation}

Both the material transport of $b$ \eqref{eq::tsw_mat_adv} and flux form transport of $B$ 
\eqref{eq::tsw_flux_adv} support skew-symmetric Hamiltonian formulations of the thermal 
shallow water equations. For the material form using prognostic variables $\boldsymbol{u}$, $h$, $b$,
we express the global energy over the domain
$\Omega$ as
\begin{equation}
	\mathcal{H}_m^{tsw}(\boldsymbol{u},h,b) = \int\Bigg(\frac{1}{2}h\boldsymbol{u}\cdot\boldsymbol{u} + \frac{1}{2}h^2b\Bigg)\mathrm{d}\Omega,
\end{equation}
for which the vector of variational derivatives $\nabla\mathcal{H}_m^{tsw}$ is expressed as
\begin{equation}
	\nabla\mathcal{H}_m^{tsw} =
	\begin{bmatrix}
		\frac{\delta\mathcal{H}_m^{tsw}}{\delta\boldsymbol{u}}\\
		\frac{\delta\mathcal{H}_m^{tsw}}{\delta h}\\
		\frac{\delta\mathcal{H}_m^{tsw}}{\delta b}
	\end{bmatrix} =
	\begin{bmatrix}
		h\boldsymbol{u}\\
		\frac{1}{2}\boldsymbol{u}\cdot\boldsymbol{u} + hb\\
		\frac{1}{2}h^2
	\end{bmatrix} = 
	\begin{bmatrix}
		\boldsymbol{F}\\
		\Phi^m\\
		T^m
	\end{bmatrix}.
\end{equation}
This leads to a skew-symmetric representation of the thermal shallow water equations as
\begin{equation}
	\frac{\partial}{\partial t}
	\begin{bmatrix}
		\boldsymbol{u}\\ h\\ b
	\end{bmatrix} = -
	\begin{bmatrix}
		q\times & \nabla & -\frac{\nabla b}{h}\\
		\nabla\cdot & 0 & 0\\
		\frac{\nabla b}{h} & 0 & 0
	\end{bmatrix} 
	\begin{bmatrix}
		\boldsymbol{F}\\
		\Phi^m\\
		T^m
	\end{bmatrix}\label{eq::tsw_ss_m}.
\end{equation}

Conversely, for the flux form, using prognostic variables $\boldsymbol{u}$, $h$, $B$, we have the global energy as
\begin{equation}
	\mathcal{H}_f^{tsw}(\boldsymbol{u},h,B) = 
	\int\Bigg(\frac{1}{2}h\boldsymbol{u}\cdot\boldsymbol{u} + \frac{1}{2}hB\Bigg)\mathrm{d}\Omega
\end{equation}
with the corresponding variational derivatives given as
\begin{equation}
	\nabla\mathcal{H}_f^{tsw} = 
	\begin{bmatrix}
		\frac{\delta\mathcal{H}_f^{tsw}}{\delta\boldsymbol{u}}\\
		\frac{\delta\mathcal{H}_f^{tsw}}{\delta h}\\
		\frac{\delta\mathcal{H}_f^{tsw}}{\delta B}
	\end{bmatrix} =
	\begin{bmatrix}
		h\boldsymbol{u}\\
		\frac{1}{2}\boldsymbol{u}\cdot\boldsymbol{u} + \frac{1}{2}B\\
		\frac{1}{2}h
	\end{bmatrix} = 
	\begin{bmatrix}
		\boldsymbol{F}\\
		\Phi^f\\
		T^f
	\end{bmatrix}.
\end{equation}
For which the skew-symmetric form of the thermal shallow water equations is expressed as
\begin{equation}
	\frac{\partial}{\partial t}
	\begin{bmatrix}
		\boldsymbol{u}\\ h\\ B
	\end{bmatrix} = -
	\begin{bmatrix}
		q\times & \nabla & \frac{B}{h}\nabla\\
		\nabla\cdot & 0 & 0\\
		\nabla\cdot(\frac{B}{h}\cdot) & 0 & 0
	\end{bmatrix} 
	\begin{bmatrix}
		\boldsymbol{F}\\
		\Phi^f\\
		T^f
	\end{bmatrix}\label{eq::tsw_ss_f}.
\end{equation}

For both \eqref{eq::tsw_ss_m} and \eqref{eq::tsw_ss_f} energy conservation arises from the skew-symmetric
nature of the right hand side operator (assuming periodic boundary conditions) and the chain rule for
the temporal derivative.
In either case multiplying both sides by $(\nabla\mathcal{H}^{tsw})^{\top}=[\boldsymbol{F}\quad\Phi\quad T]$, 
and integrating over the domain $\Omega$ leads to the expression
\begin{equation}\int
	\begin{bmatrix}
		\frac{\delta\mathcal{H}^{tsw}}{\delta\boldsymbol{u}} &
		\frac{\delta\mathcal{H}^{tsw}}{\delta h} &
		\frac{\delta\mathcal{H}^{tsw}}{\delta a}
	\end{bmatrix}
	\begin{bmatrix}
		\frac{\partial\boldsymbol{u}}{\partial t} \\
		\frac{\partial h}{\partial t} \\
		\frac{\partial a}{\partial t}
	\end{bmatrix}\mathrm{d}\Omega = 
	\int\Bigg(\frac{\delta\mathcal{H}^{tsw}}{\delta\boldsymbol{u}}\cdot\frac{\partial\boldsymbol{u}}{\partial t} +
		\frac{\delta\mathcal{H}^{tsw}}{\delta h}\cdot\frac{\partial h}{\partial t} +
		\frac{\delta\mathcal{H}^{tsw}}{\delta a}\cdot\frac{\partial a}{\partial t}\Bigg)\mathrm{d}\Omega
	=\frac{\mathrm{d}\mathcal{H}^{tsw}}{\mathrm{d}t} = 0
\end{equation}
for $a\in \{b,B\}$.

Additionally, both skew-symmetric operators have a null-space, which contains the 
variational derivatives of the additional invariants of the system. These include
the mass, $\mathcal{M}^{tsw}=\int h\mathrm{d}\Omega$, the vorticity, $\mathcal{W}^{tsw}=\int (hq-f)\mathrm{d}\Omega$,
the density weighted buoyancy, $\mathcal{B}^{tsw} = \int hb\mathrm{d}\Omega$ \cite{Eldred19}, as well as
the higher order moments of the buoyancy \cite{Dellar03} such as the
entropy $\mathcal{S}^{tsw}=\frac{1}{2}\int hb^2\mathrm{d}\Omega = \int \frac{B^2}{2h}\mathrm{d}\Omega$. 
While the preservation of all of these quantities is dynamically important, we are 
principally concerned in the present study with the conservation of entropy, since 
this is known to have eigenvalues which are all $\ge 0$, and would therefore ensure the stability 
of thermal processes.

The variational derivatives of the entropy for the material form are given as
$\nabla\mathcal{S}_m^{tsw} = [\boldsymbol{0}, \frac{1}{2}b^2, hb]^{\top}$, such 
that this is conserved as
	\begin{align}
		\frac{\mathrm{d}\mathcal{S}_m^{tsw}}{\mathrm{d}t} &= 
		\int
	\begin{bmatrix}
		\boldsymbol{0} & \frac{1}{2}b^2 & hb
	\end{bmatrix}
	\begin{bmatrix}
		\frac{\partial\boldsymbol{u}}{\partial t} \\
		\frac{\partial h}{\partial t} \\
		\frac{\partial a}{\partial t}
	\end{bmatrix} \mathrm{d}\Omega \notag \\
		&= -\int
	\begin{bmatrix}
		\boldsymbol{0} & \frac{1}{2}b^2 & hb
	\end{bmatrix}
	\begin{bmatrix}
		q\times & \nabla & -\frac{\nabla b}{h}\\
		\nabla\cdot & 0 & 0\\
		\frac{\nabla b}{h} & 0 & 0
	\end{bmatrix} 
	\begin{bmatrix}
		\boldsymbol{F}\\
		\Phi^m\\
		T^m
	\end{bmatrix} \mathrm{d}\Omega \notag \\
		&= \int
	\begin{bmatrix}
		\frac{1}{2}\nabla b^2 - b\nabla b & 0 & 0
	\end{bmatrix}
	\begin{bmatrix}
		\boldsymbol{F}\\
		\Phi^m\\
		T^m
	\end{bmatrix}\mathrm{d}\Omega \notag \\
		&= 
		\int\Bigg(\frac{1}{2}\nabla(b^2)\cdot\boldsymbol{F} - b\nabla(b)\cdot\boldsymbol{F}\Bigg)\mathrm{d}\Omega
		=0,\label{eq::entropy_cons_m}
	\end{align}
where we have applied integration by parts to the first term in the row vector (assuming periodic boundary conditions),
and then invoked the chain rule in order to show conservation.
Similarly for the flux form, the variational derivatives for the corresponding 
entropy are given as $\nabla\mathcal{S}_f^{tsw} = [\boldsymbol{0}, \frac{-B^2}{2h^2}, \frac{B}{h}]^{\top}$, for which 
the entropy is conserved as 
\begin{align}
	\frac{\mathrm{d}\mathcal{S}_f^{tsw}}{\mathrm{d}t} &= 
	\int
	\begin{bmatrix}
		\boldsymbol{0} & -\frac{B^2}{2h^2} & \frac{B}{h}
	\end{bmatrix}
	\begin{bmatrix}
		\frac{\partial\boldsymbol{u}}{\partial t} \\
		\frac{\partial h}{\partial t} \\
		\frac{\partial a}{\partial t}
	\end{bmatrix} \mathrm{d}\Omega \notag \\
	&= -\int
	\begin{bmatrix}
		\boldsymbol{0} & -\frac{B^2}{2h^2} & \frac{B}{h}
	\end{bmatrix}
	\begin{bmatrix}
		q\times & \nabla & \frac{B}{h}\nabla\\
		\nabla\cdot & 0 & 0\\
		\nabla\cdot(\frac{B}{h}\cdot) & 0 & 0
	\end{bmatrix} 
	\begin{bmatrix}
		\boldsymbol{F}\\
		\Phi^f\\
		T^f
	\end{bmatrix}\mathrm{d}\Omega \notag \\
	&= \int
	\begin{bmatrix}
		-\frac{1}{2}\nabla\Big(\frac{B^2}{h^2}\Big) + \frac{B}{h}\nabla\Big(\frac{B}{h}\Big) & 0 & 0
	\end{bmatrix}
	\begin{bmatrix}
		\boldsymbol{F}\\
		\Phi^f\\
		T^f
	\end{bmatrix}\mathrm{d}\Omega \notag \\
	&=
	\int\Bigg(-\frac{1}{2}\nabla\Big(\frac{B^2}{h^2}\Big)\cdot\boldsymbol{F} + 
	\frac{B}{h}\nabla\Big(\frac{B}{h}\Big)\cdot\boldsymbol{F}\Bigg)\mathrm{d}\Omega
	= 0,\label{eq::entropy_cons_f}
\end{align}
where above we have applied integration by parts to both terms in the row vector and again
invoked the chain rule to establish conservation.
In both cases, entropy conservation is derived from the chain rule. This is problematic 
in the discrete form if $h$ and $b$ or $B$ are represented by discontinuous function spaces, 
for which the chain rule is not defined at element boundaries, or if inexact integration
is used, for which the chain rule also fails. 

In order to negate this issue and ensure that entropy is conserved independent of the chain
rule, we may alternatively construct our solution as an average of \eqref{eq::tsw_ss_m} and 
\eqref{eq::tsw_ss_f} as
\begin{equation}\label{eq::tsw_mixed}
	\frac{\partial}{\partial t}
	\begin{bmatrix}\boldsymbol{u} \\ h \\ \frac{B}{2} \\ \frac{b}{2} \end{bmatrix}=
		-\frac{1}{2}\begin{bmatrix}q\times & \nabla & b\nabla & \boldsymbol{0} \\
			\nabla\cdot & 0 & 0 & 0 \\
			\nabla\cdot(b\cdot) & 0 & 0 & 0 \\
			0 & 0 & 0 & 0\end{bmatrix}
			\begin{bmatrix}\boldsymbol{F} \\ \Phi^f \\ T^f \\ 0 \end{bmatrix}
		-\frac{1}{2}\begin{bmatrix}q\times & \nabla & \boldsymbol{0} & -\frac{1}{h}\nabla\Big(\frac{B}{h}\Big) \\
			\nabla\cdot & 0 & 0 & 0 \\
			0 & 0 & 0 & 0 \\
			\frac{1}{h}\nabla\Big(\frac{B}{h}\Big) & 0 & 0 & 0\end{bmatrix}
			\begin{bmatrix}\boldsymbol{F} \\ \Phi^m \\ 0 \\ T^m \end{bmatrix},
\end{equation}
where we solve for both the buoyancy and density weighted buoyancy prognostically as
an alternative to diagnosing one from the other.
Note that the above formulation differs from a weighted average of \eqref{eq::tsw_ss_m} and
\eqref{eq::tsw_ss_f} in one respect: the contents of the skew-symmetric matrix for the material 
form, \eqref{eq::tsw_ss_m} are instead of $b$ using $\frac{B}{h}$, while the skew-symmetric 
operator in \eqref{eq::tsw_ss_f} are instead of $B$ using $hb$. This will become important in
the discrete form, since it will allow for a tighter coupling of the material and flux form
transport equations.
Now pre-multiplying by $\frac{1}{2}(\nabla\mathcal{S}_m^{tsw} + \nabla\mathcal{S}_f^{tsw})^{\top}$
and recalling that $b=\frac{B}{h}$ leads to 
cancellations separately between the first terms in \eqref{eq::entropy_cons_m} and \eqref{eq::entropy_cons_f}
and the second terms in \eqref{eq::entropy_cons_m} and \eqref{eq::entropy_cons_f} such
that $\mathcal{S}_m^{tsw} = \mathcal{S}_f^{tsw}$ is conserved independently of the chain rule. 

More succinctly, we may combine \eqref{eq::tsw_cont} and \eqref{eq::tsw_mat_adv} to eliminate
the material advection equation for $b$ from the above system to give
\begin{equation}\label{eq::tsw_coupled}
	\frac{\partial}{\partial t}
	\begin{bmatrix}\boldsymbol{u} \\ h \\ B \end{bmatrix}=
		-\begin{bmatrix}q\times & \nabla & \frac{1}{2}b\nabla + \frac{1}{2}\nabla(b\cdot) - \frac{1}{2}\nabla b \\
		\nabla\cdot & 0 & 0 \\
			\frac{1}{2}\nabla\cdot(b\cdot) + \frac{1}{2}b\nabla\cdot + \frac{1}{2}\nabla b & 0 & 0 \\
		\end{bmatrix}
		\begin{bmatrix}\boldsymbol{F} \\ \Phi^f \\ T^f \end{bmatrix}
\end{equation}
as has been previously done for an entropy conserving discretisation of the thermal shallow
water equations using a discontinuous Galerkin method \cite{Ricardo23b}. Once again, pre-multiplying
\eqref{eq::tsw_coupled} by $(\nabla\mathcal{H}_f^{tsw})^{\top}$ allows for the conservation of energy
due to the skew-symmetry of the right hand side and integration by parts as
\begin{align}
	\frac{\mathrm{d}\mathcal{H}_f^{tsw}}{\mathrm{d}t} = -
	\int\begin{bmatrix}\boldsymbol{F} & \Phi^f & T^f \end{bmatrix} 
		\begin{bmatrix}q\times & \nabla & \frac{1}{2}b\nabla + \frac{1}{2}\nabla(b\cdot) - \frac{1}{2}\nabla b \\
		\nabla\cdot & 0 & 0 \\
			\frac{1}{2}\nabla\cdot(b\cdot) + \frac{1}{2}b\nabla\cdot + \frac{1}{2}\nabla b & 0 & 0 \\
		\end{bmatrix}
		\begin{bmatrix}\boldsymbol{F} \\ \Phi^f \\ T^f \end{bmatrix}\mathrm{d}\Omega &= \notag \\
			-\int\Bigg(\boldsymbol{F}q\times\boldsymbol{F} + \boldsymbol{F}\nabla\Phi^f + 
	%\frac{1}{2}\Big(b\boldsymbol{F}\cdot\nabla T^f + T^f\nabla\cdot(b\boldsymbol{F}) - T^f\boldsymbol{F}\cdot\nabla b\Big) &\notag\\
	\frac{1}{2}\Big(b\boldsymbol{F}\cdot\nabla T^f + \blue{\boldsymbol{F}\cdot\nabla(bT^f) } - T^f\boldsymbol{F}\cdot\nabla b\Big) &\notag\\
	+\Phi^f\nabla\cdot\boldsymbol{F} + \frac{1}{2}\Big(
	T^f\nabla\cdot(b\boldsymbol{F}) + b\blue{T^f}\nabla\cdot\boldsymbol{F} + T^f\boldsymbol{F}\cdot\nabla b\Big)\Bigg)
	\mathrm{d}\Omega &= 0,
\end{align}
while pre-multiplying by $(\nabla\mathcal{S}_f^{tsw})^{\top}$ allows for the conservation of entropy
via integration by parts but without necessitating the invocation of the chain rule as
\begin{align}
	\frac{\mathrm{d}\mathcal{S}_f^{tsw}}{\mathrm{d}t} = -
	\int
	\begin{bmatrix}\boldsymbol{0} & -\frac{1}{2}b^2 & b \end{bmatrix} 
		\begin{bmatrix}q\times & \nabla & \frac{1}{2}b\nabla + \frac{1}{2}\nabla(b\cdot) - \frac{1}{2}\nabla b \\
		\nabla\cdot & 0 & 0 \\
			\frac{1}{2}\nabla\cdot(b\cdot) + \frac{1}{2}b\nabla\cdot + \frac{1}{2}\nabla b & 0 & 0 \\
		\end{bmatrix}
		\begin{bmatrix}\boldsymbol{F} \\ \Phi^f \\ T^f \end{bmatrix}\mathrm{d}\Omega &= \notag \\
			-\frac{1}{2}\int\Bigg(-b^2\nabla\cdot\boldsymbol{F} + b\nabla\cdot(b\boldsymbol{F}) +
			b^2\nabla\cdot\boldsymbol{F} + b\boldsymbol{F}\cdot\nabla b
	\Bigg)\mathrm{d}\Omega &= 0.
\end{align}

\subsection{3D compressible Euler equations}

Extending the above formulation to the 3D compressible Euler equations is remarkably simple. This is
because the 3D compressible Euler equations have the same skew-symmetric structure as the thermal
shallow water equations as given in \eqref{eq::tsw_ss_m} and \eqref{eq::tsw_ss_f}. The only thing that
differs is the energy and is variational derivatives, which for the flux form 
using prognostic variables $\boldsymbol{u}$, $\rho$, $\Theta$, is given as 
\begin{equation}
	\mathcal{H}_f^{ce}(\boldsymbol{u},\rho,\Theta) = \int\Bigg(\frac{1}{2}\rho\boldsymbol{u}\cdot\boldsymbol{u} + 
	\rho gz + \frac{c_v}{c_p}\Theta\Pi\Bigg)\mathrm{d}\Omega,
\end{equation}
where $g$ is gravity, $z$ is height, $\rho$ is the density, $\Theta$ the density weighted weighted
potential temperature (and $\theta:=\frac{\Theta}{\rho}$ is the potential temperature) and $c_v$, $c_p$ are
the specific heats at constant volume and pressure respectively. The Exner pressure is denoted by
$\Pi$, which is derived from the ideal gas law as
\begin{equation}\label{eq::eos}
	\Pi = c_p\Bigg(\frac{R\Theta}{p_0}\Bigg)^{R/c_v}, 
\end{equation}
with $p_0$ being the reference pressure at the surface and $R=c_p-c_v$ is the ideal gas constant.

The variational derivatives of the energy for the 3D compressible Euler equations with the flux form
transport of density weighted potential temperature are then given as
\begin{equation}\label{eq::ce_grad_H}
	\nabla\mathcal{H}_f^{ce} = 
	\begin{bmatrix}\rho\boldsymbol{u} \\ \frac{1}{2}\boldsymbol{u}\cdot\boldsymbol{u} + gz \\
	c_p\Big(\frac{R\Theta}{p_0}\Big)^{\frac{R}{c_v}}
	\end{bmatrix} = 
	\begin{bmatrix}
		\boldsymbol{U} \\ \Phi \\ \Pi
	\end{bmatrix}.
\end{equation}
Note that we have omitted superscripts for the equation set and subscripts for flux form from the
above definitions for the sake of brevity. The 3D compressible Euler equations can then be expressed 
in a direct analogue of \eqref{eq::tsw_ss_f} as
\begin{equation}
	\frac{\partial}{\partial t}
	\begin{bmatrix}
		\boldsymbol{u}\\ \rho\\ \Theta
	\end{bmatrix} = -
	\begin{bmatrix}
		q\times & \nabla & \theta\nabla\\
		\nabla\cdot & 0 & 0\\
		\nabla\cdot(\theta\cdot) & 0 & 0
	\end{bmatrix} 
	\begin{bmatrix}
		\boldsymbol{U}\\
		\Phi\\
		\Pi
	\end{bmatrix}\label{eq::ce_ss}.
\end{equation}
Alternatively we may then construct an analogue of the coupled thermal shallow water system 
\eqref{eq::tsw_coupled} for the 3D compressible Euler equations as
\begin{equation}\label{eq::ce_coupled}
	\frac{\partial}{\partial t}
	\begin{bmatrix}\boldsymbol{u} \\ \rho \\ \Theta \end{bmatrix}=
		-\begin{bmatrix}\boldsymbol{q}\times & \nabla & \frac{1}{2}\theta\nabla + \frac{1}{2}\nabla(\theta\cdot) - \frac{1}{2}\nabla\theta \\
		\nabla\cdot & 0 & 0 \\
			\frac{1}{2}\nabla\cdot(\theta\cdot) + \frac{1}{2}\theta\nabla\cdot + \frac{1}{2}\nabla\theta & 0 & 0 \\
		\end{bmatrix}
		\begin{bmatrix}\boldsymbol{U} \\ \Phi \\ \Pi \end{bmatrix}.
\end{equation}
The corresponding entropy is given as 
$\mathcal{S}^{ce} = \frac{1}{2}\int\rho\theta^2\mathrm{d}\Omega$, 
and its variational derivatives as 
$\nabla\mathcal{S}^{ce}=[\boldsymbol{0}, -\frac{1}{2}\theta^2, \theta]^{\top}$. 
Once again, pre-multiplication of \eqref{eq::ce_coupled} by the transpose of this vector allows 
for the conservation of $\mathcal{S}^{ce}$ without resorting to the chain rule.

\section{Mixed finite element discretisations}

In this study we apply a mixed compatible finite element discretisation in order to preserve
the energy and entropy conserving properties of the continuous formulations discussed in the 
previous section to the discrete form. 

\subsection{Thermal shallow water equations}

For the thermal shallow water equations we discretise the two dimensional domain $\Omega^2$
using the finite dimensional subspaces $\mathbb{V}_0\subset H^1(\Omega^2)$, 
$\mathbb{V}_1\subset H(\mathrm{div},\Omega^2)$, $\mathbb{V}_2\subset L^2(\Omega^2)$. 
These subspaces are spanned by the basis functions $\psi_h\in\mathbb{V}_0$, 
$\boldsymbol{v}_h\in\mathbb{V}_1$ and $\phi_h\in\mathbb{V}_2$ respectively. The 
compatible relationships between these subspaces may then be expressed via a discrete
de Rham complex as \cite{AFW06,Cotter23}
\begin{equation}\label{eq::deRham_2d}
	\mathbb{R}\longrightarrow\mathbb{V}_0\stackrel{\nabla^{\perp}}{\longrightarrow}
	\mathbb{V}_1\stackrel{\nabla\cdot}{\longrightarrow}\mathbb{V}_2\longrightarrow 0,
\end{equation}
such that the vector identity $\nabla\cdot\nabla^{\perp}:=0$ is preserved discretely. 
The discrete analogues of the \emph{curl} and \emph{gradient} operators are then applied
as adjoints of the strong form $\nabla^{\perp}$ and $\nabla\cdot$ operators 
\cite{MC14,LPG18,Eldred19}.

In order to discretely solve the entropy conserving thermal shallow water equations
we first introduce the operators:
\begin{subequations}
	\begin{alignat}{5}
		\boldsymbol{\mathsf{M}}_0&:\mathbb{V}_0\rightarrow\mathbb{V}_0'&&\quad\mathrm{for\ which}\quad
		\boldsymbol{\mathsf{M}}_0\psi_h&&:=\langle\psi_h,\cdot\rangle,\quad&&\forall\psi_h\in\mathbb{V}_0,\\
		\boldsymbol{\mathsf{M}}_1&:\mathbb{V}_1\rightarrow\mathbb{V}_1'&&\quad\mathrm{for\ which}\quad
		\boldsymbol{\mathsf{M}}_1\boldsymbol{v}_h&&:=\langle\boldsymbol{v}_h,\cdot\rangle,\quad&&\forall\boldsymbol{v}_h\in\mathbb{V}_1,\\
		\boldsymbol{\mathsf{M}}_2&:\mathbb{V}_2\rightarrow\mathbb{V}_2'&&\quad\mathrm{for\ which}\quad
		\boldsymbol{\mathsf{M}}_2\phi_h&&:=\langle\phi_h,\cdot\rangle,\quad&&\forall\phi_h\in\mathbb{V}_2,\\
		\boldsymbol{\mathsf{D}}_2&:\mathbb{V}_1\rightarrow\mathbb{V}_2'&&\quad\mathrm{for\ which}\quad
		\boldsymbol{\mathsf{D}}_2\boldsymbol{v}_h&&:=\langle\nabla\cdot\boldsymbol{v}_h,\cdot\rangle,\quad&&\forall\boldsymbol{v}_h\in\mathbb{V}_1,\\
		\boldsymbol{\mathsf{R}}_1&:\mathbb{V}_0\rightarrow\mathbb{V}_1'&&\quad\mathrm{for\ which}\quad
		\boldsymbol{\mathsf{R}}_1\psi_h&&:=\langle\nabla^{\perp}\psi_h,\cdot\rangle,\quad&&\forall\psi_h\in\mathbb{V}_0,\\
		\boldsymbol{\mathsf{C}}_1&:\mathbb{V}_0\otimes\mathbb{V}_1\rightarrow\mathbb{V}_1'&&\quad\mathrm{for\ which}\quad
		\boldsymbol{\mathsf{C}}_1(\psi_h,\boldsymbol{v}_h)&&:=\langle\psi_h\times\boldsymbol{v}_h,\cdot\rangle,\quad&&\forall\psi_h\in\mathbb{V}_0,\boldsymbol{v}_h\in\mathbb{V}_1,\\
		\boldsymbol{\mathsf{M}}_{0*}&:\mathbb{V}_0\otimes\mathbb{V}_2\rightarrow\mathbb{V}_0'&&\quad\mathrm{for\ which}\quad
		\boldsymbol{\mathsf{M}}_{0*}(\psi_h,\phi_h)&&:=\langle\phi_h\psi_h,\cdot\rangle,\quad&&\forall\psi_h\in\mathbb{V}_0,\phi_h\in\mathbb{V}_2,\\
		\boldsymbol{\mathsf{M}}_{1*}&:\mathbb{V}_1\otimes\mathbb{V}_2\rightarrow\mathbb{V}_1'&&\quad\mathrm{for\ which}\quad
		\boldsymbol{\mathsf{M}}_{1*}(\boldsymbol{v}_h,\phi_h)&&:=\langle\phi_h\boldsymbol{v}_h,\cdot\rangle,\quad&&\forall\boldsymbol{v}_h\in\mathbb{V}_1,\phi_h\in\mathbb{V}_2,\\
		\boldsymbol{\mathsf{M}}_{2*}&:\mathbb{V}_2\otimes\mathbb{V}_2\rightarrow\mathbb{V}_2'&&\quad\mathrm{for\ which}\quad
		\boldsymbol{\mathsf{M}}_{2*}(\phi_h,\phi_k)&&:=\langle\phi_h\phi_k,\cdot\rangle,\quad&&\forall\phi_h,\phi_k\in\mathbb{V}_2,\\
		\boldsymbol{\mathsf{K}}_2&:\mathbb{V}_1\otimes\mathbb{V}_1\rightarrow\mathbb{V}_2'&&\quad\mathrm{for\ which}\quad
		\boldsymbol{\mathsf{K}}_2(\boldsymbol{v}_h,\boldsymbol{v}_k)&&:=\langle\boldsymbol{v}_h\cdot\boldsymbol{v}_k,\cdot\rangle,\quad&&\forall\boldsymbol{v}_h,\boldsymbol{v}_k\in\mathbb{V}_1,
	\end{alignat}
\end{subequations}
where primes denote corresponding algebraic dual spaces such that the operators above map trial functions into dual 
test spaces, and angled brackets denote inner products as $\langle a,b\rangle=\int ab\mathrm{d}\Omega^2$, and likewise for 
tri-linear operators. For the matrices above that represent differential operators, $\boldsymbol{\mathsf{R}}_1$, 
$\boldsymbol{\mathsf{D}}_2$, these are applied to the trial function, such that they are \emph{strong form} operators as detailed 
in the de-Rham complex \eqref{eq::deRham_2d}. Their \emph{weak form} adjoints, for which the differential operator is applied 
to the test function, are given as the transpose of the associated strong form operator.

\subsubsection{Thermal shallow water equations: mixed formulation}

Having defined the discrete function spaces and operators, we now proceed to describe the semi-discrete 
analogue of the mixed material-flux form of the thermal shallow water equations \eqref{eq::tsw_mixed}
for the discrete solution variables $\boldsymbol{u}_h\in\mathbb{V}_1$, $h_h, B_h, b_h\in\mathbb{V}_2$.
Firstly the diagnostic equations for the discrete potential vorticity, $q_h\in\mathbb{V}_0$ and 
variational derivatives $\boldsymbol{F}_h\in\mathbb{V}_1$, 
$\Phi^m_h, \Phi^f_h, T^m_h, T^f_h\in\mathbb{V}_2$ are given as
\begin{subequations}
	\begin{align}
		\boldsymbol{\mathsf{M}}_{0*}(q_h,h_h) &= -\boldsymbol{\mathsf{R}}_1^{\top}\boldsymbol{u}_h + \boldsymbol{\mathsf{M}}_0f_h \label{eq::q_diag}\\
		\boldsymbol{\mathsf{M}}_1\boldsymbol{F}_h &= \boldsymbol{\mathsf{M}}_{1*}(\boldsymbol{u}_h,h_h) \label{eq::F_diag}\\
		\boldsymbol{\mathsf{M}}_2\Phi^m_h &= \frac{1}{2}\boldsymbol{\mathsf{K}}_2(\boldsymbol{u}_h,\boldsymbol{u}_h) + 
		\boldsymbol{\mathsf{M}}_{2*}(h_h,b_h) \\
		\boldsymbol{\mathsf{M}}_2\Phi^f_h &= \frac{1}{2}\boldsymbol{\mathsf{K}}_2(\boldsymbol{u}_h,\boldsymbol{u}_h) + 
		\frac{1}{2}\boldsymbol{\mathsf{M}}_2B_h \\
		\boldsymbol{\mathsf{M}}_2T^m_h &= \frac{1}{2}\boldsymbol{\mathsf{M}}_{2*}(h_h,h_h) \\
		\boldsymbol{\mathsf{M}}_2T^f_h &= \frac{1}{2}\boldsymbol{\mathsf{M}}_2h_h,
	\end{align}
\end{subequations}
where $f_h\in\mathbb{V}_0$ is the discrete representation of the Coriolis term.
We also compute a discrete analogue of the buoyancy gradient per unit depth, $\boldsymbol{G}:= \frac{1}{h}\nabla b'$,
where $b' = \frac{B}{h}$, as appears in the skew symmetric operator for the material form of the thermal shallow water equations 
\eqref{eq::tsw_ss_m} for $\boldsymbol{G}_h\in\mathbb{V}_1$, $b_h'\in\mathbb{V}_2$ as
\begin{subequations}\label{eq::G_from_B}
	\begin{align}
		\boldsymbol{\mathsf{M}}_{2*}(h_h,b_h') &= \boldsymbol{\mathsf{M}}_2B_h \label{eq::b_prime}\\
		\boldsymbol{\mathsf{M}}_{1*}(\boldsymbol{G}_h,h_h) &= -\boldsymbol{\mathsf{D}}_2^{\top}b_h'.\label{eq::G_diag}
	\end{align}
\end{subequations}

Once we have computed the necessary diagnostic quantities, we may solve for the prognostic variables in
the semi-discrete form of \eqref{eq::tsw_mixed} as
\begin{subequations}
	\begin{align}
		\boldsymbol{\mathsf{M}}_1\frac{\mathrm{d}\boldsymbol{u}_h}{\mathrm{d} t} &=
		-\boldsymbol{\mathsf{C}}_1(q_h,\boldsymbol{F}_h) + \frac{1}{2}\boldsymbol{\mathsf{D}}_2^{\top}(\Phi^m_h + \Phi^f_h) 
		+\frac{1}{2}\boldsymbol{\mathsf{K}}_2^{\top}(\boldsymbol{G}_h,T_h^m) 
		+\frac{1}{2}\boldsymbol{\mathsf{M}}_{1*}(\boldsymbol{\mathsf{M}}_1^{-1}\boldsymbol{\mathsf{D}}_2^{\top}T_h^f,b_h) \label{eq::tsw_mixed_disc_u}\\
		\boldsymbol{\mathsf{M}}_2\frac{\mathrm{d} h_h}{\mathrm{d} t} &= -\boldsymbol{\mathsf{D}}_2\boldsymbol{F}_h \label{eq::tsw_mixed_disc_h}\\
		\boldsymbol{\mathsf{M}}_2\frac{\mathrm{d} b_h}{\mathrm{d} t} &= -\boldsymbol{\mathsf{K}}_2(\boldsymbol{G}_h,\boldsymbol{F}_h) \label{eq::tsw_mixed_disc_b}\\
		\boldsymbol{\mathsf{M}}_2\frac{\mathrm{d} B_h}{\mathrm{d} t} &= 
		-\boldsymbol{\mathsf{D}}_2\boldsymbol{\mathsf{M}}_1^{-1}\boldsymbol{\mathsf{M}}_{1*}(\boldsymbol{F}_h,b_h).\label{eq::tsw_mixed_disc_B}
	\end{align}\label{eq::tsw_mixed_disc}
\end{subequations}
Note the presence of the $\mathbb{V}_1$ mass matrix inverse $\boldsymbol{\mathsf{M}}_1^{-1}$ in \eqref{eq::tsw_mixed_disc_u}
and \eqref{eq::tsw_mixed_disc_B} for the computation of the weak gradient of $T_h^f$ and the buoyancy flux respectively. This implies 
a projection onto the $\mathbb{V}_1$ space via an implicit solve, which adds to the expense of the algorithm. Additional projections
onto the $\mathbb{V}_0$ and $\mathbb{V}_1$ spaces are required for the diagnostic equations in \eqref{eq::q_diag}, \eqref{eq::F_diag} 
and \eqref{eq::G_diag} and the prognostic equation in \eqref{eq::tsw_mixed_disc_u}. While there are additional projections for the 
various diagnostic and prognostic equations onto the $\mathbb{V}_2$ space, these are trivial since the basis functions for this
space are discontinuous in all dimensions, and so these may be computed element-wise.

Left-multiplication of \eqref{eq::tsw_mixed_disc} by the vector 
$[\boldsymbol{F}_h^{\top}, \frac{1}{2}(\Phi_h^{m\top} + \Phi_h^{f\top}), \frac{1}{2}T_h^{m\top}, \frac{1}{2}T_h^{f\top}]$
leads to the cancellation of all right hand side terms as
\begin{subequations}
	\begin{align}
		-\boldsymbol{F}_h^{\top}\boldsymbol{\mathsf{C}}_1(q_h,\boldsymbol{F}_h) + \frac{1}{2}\boldsymbol{F}_h^{\top}\boldsymbol{\mathsf{D}}_2^{\top}(\Phi^m_h + \Phi^f_h) 
		+\frac{1}{2}\boldsymbol{F}_h^{\top}\boldsymbol{\mathsf{K}}_2^{\top}(\boldsymbol{G}_h,T_h^m) 
		+\frac{1}{2}\boldsymbol{F}_h^{\top}\boldsymbol{\mathsf{M}}_{1*}(\boldsymbol{\mathsf{M}}_1^{-1}\boldsymbol{\mathsf{D}}_2^{\top}T_h^f,b_h) &\notag \\
		-\frac{1}{2}(\Phi_h^{m\top} + \Phi_h^{f\top})\boldsymbol{\mathsf{D}}_2\boldsymbol{F}_h
		-\frac{1}{2}T_h^{m\top}\boldsymbol{\mathsf{K}}_2(\boldsymbol{G}_h,\boldsymbol{F}_h)
		-\frac{1}{2}T_h^{f\top}\boldsymbol{\mathsf{D}}_2\boldsymbol{\mathsf{M}}_1^{-1}\boldsymbol{\mathsf{M}}_{1*}(\boldsymbol{F}_h,b_h) &= \\
		-\langle\boldsymbol{F}_h,\boldsymbol{q}_h\times\boldsymbol{F}_h\rangle + \frac{1}{2}\langle\nabla\cdot\boldsymbol{F}_h,\Phi_h^m+\Phi_h^f\rangle + 
		\frac{1}{2}\langle\boldsymbol{F}_h\cdot\boldsymbol{G}_h,T_h^m\rangle + \frac{1}{2}\langle\nabla\cdot(b_h\boldsymbol{F}_h),T_h^f\rangle &\notag \\
		-\frac{1}{2}\langle\Phi_h^m+\Phi_h^f,\nabla\cdot\boldsymbol{F}_h\rangle - \frac{1}{2}\langle T_h^m,\boldsymbol{G}_h\cdot\boldsymbol{F}_h\rangle 
		-\frac{1}{2}\langle T_h^f,\nabla\cdot(b_h\boldsymbol{F}_h)\rangle &= 0,
	\end{align}
\end{subequations}
and results in the expression
\begin{equation}\label{eq::temporal_en_cons}
	\frac{1}{2}\Bigg(\boldsymbol{F}_h^{\top}\boldsymbol{\mathsf{M}}_1\frac{\mathrm{d}\boldsymbol{u}_h}{\mathrm{d} t} + 
	\Phi_h^{m\top}\boldsymbol{\mathsf{M}}_2\frac{\mathrm{d} h_h}{\mathrm{d} t} + 
	T_h^{m\top}\boldsymbol{\mathsf{M}}_2\frac{\mathrm{d} b_h}{\mathrm{d} t} + 
	\boldsymbol{F}_h^\top\boldsymbol{\mathsf{M}}_1\frac{\mathrm{d}\boldsymbol{u}_h}{\mathrm{d} t} + 
	\Phi_h^{f\top}\boldsymbol{\mathsf{M}}_2\frac{\mathrm{d} h_h}{\mathrm{d} t} + 
	T_h^{f\top}\boldsymbol{\mathsf{M}}_2\frac{\mathrm{d} B_h}{\mathrm{d} t}\Bigg) = 
	\frac{1}{2}\Bigg(\frac{\mathrm{d}\mathcal{H}_{m,h}^{tsw}}{\mathrm{d} t} + \frac{\mathrm{d}\mathcal{H}_{f,h}^{tsw}}{\mathrm{d} t}\Bigg) = 0.
\end{equation}
In order to further conserve energy in time as well as space, the functional derivatives on the left
hand side in the above expression should be exactly integrated in time up to the order of the temporal 
discretisation \cite{BC18,Eldred19,Lee21,LP21}.

The discrete functional derivatives of the entropy are expressed as
\begin{subequations}
	\begin{align}
		\nabla\mathcal{S}_{m,h}^{tsw} &= 
		\begin{bmatrix} \boldsymbol{0} \\ 
			\frac{1}{2}\Big(\boldsymbol{\mathsf{M}}_{2*}(b_h,b_h)\Big)^{\top}\boldsymbol{\mathsf{M}}_2^{-1} \\
			\Big(\boldsymbol{\mathsf{M}}_{2*}(h_h,b_h)\Big)^{\top}\boldsymbol{\mathsf{M}}_2^{-1} 
		\end{bmatrix} \label{eq::grad_S_m_disc}\\
		\nabla\mathcal{S}_{f,h}^{tsw} &=
		\begin{bmatrix} \boldsymbol{0} \\ 
			-\frac{1}{2}\Big(\boldsymbol{\mathsf{M}}_{2*}(B_h,B_h)\Big)^{\top}
			\boldsymbol{\mathsf{M}}_{2*}^{-1}(h_h,\boldsymbol{\mathsf{M}}_2\boldsymbol{\mathsf{M}}_{2*}^{-1}(h_h,\cdot)) \\
			\Big(\boldsymbol{\mathsf{M}}_2B_h\Big)^{\top}\boldsymbol{\mathsf{M}}_{2*}^{-1}(h_h,\cdot)
		\end{bmatrix}.\label{eq::grad_S_f_disc}
	\end{align}
\end{subequations}
Left multiplication of \eqref{eq::tsw_mixed_disc_h} by an average of the second rows in \eqref{eq::grad_S_m_disc} 
and \eqref{eq::grad_S_f_disc} results in a temporal change of entropy (i.e. a loss of conservation) as
\begin{equation}\label{eq::tsw_mixed_S1}
	\frac{1}{2}\Big(\boldsymbol{\mathsf{M}}_{2*}(b_h,b_h)\Big)^{\top}\boldsymbol{\mathsf{M}}_2^{-1}\boldsymbol{\mathsf{D}}_2\boldsymbol{F}_h -
	\frac{1}{2}\Big(\boldsymbol{\mathsf{M}}_{2*}(B_h,B_h)\Big)^{\top}
	\boldsymbol{\mathsf{M}}_{2*}^{-1}(h_h,\boldsymbol{\mathsf{M}}_2\boldsymbol{\mathsf{M}}_{2*}^{-1}(h_h,\boldsymbol{\mathsf{D}}_2\boldsymbol{F}_h)) = 
	\frac{1}{2}\Big\langle b_h^2,\nabla\cdot\boldsymbol{F}_h\Big\rangle - 
	\frac{1}{2}\Big\langle\frac{B_h^2}{h_h^2},\nabla\cdot\boldsymbol{F}_h\Big\rangle.
\end{equation}
Similarly, left multiplication of \eqref{eq::tsw_mixed_disc_b} by the bottom 
row in \eqref{eq::grad_S_m_disc} and \eqref{eq::tsw_mixed_disc_B} by the bottom row in \eqref{eq::grad_S_f_disc} and adding then 
gives
\begin{align}\label{eq::tsw_mixed_S2}
	-\Big(\boldsymbol{\mathsf{M}}_{2*}(h_h,b_h)\Big)^{\top}\boldsymbol{\mathsf{M}}_2^{-1} 
	\boldsymbol{\mathsf{K}}_2(\boldsymbol{F}_h,\boldsymbol{\mathsf{M}}_{1*}^{-1}(
	\boldsymbol{\mathsf{D}}_2^{\top}\boldsymbol{\mathsf{M}}_{2*}^{-1}(h_h,\boldsymbol{\mathsf{M}}_2B_h),h_h))
	+\Big(\boldsymbol{\mathsf{M}}_2B_h\Big)^{\top}\boldsymbol{\mathsf{M}}_{2*}^{-1}(h_h,
	\boldsymbol{\mathsf{D}}_2\boldsymbol{\mathsf{M}}_1^{-1}\boldsymbol{\mathsf{M}}_{1*}(\boldsymbol{F}_h,b_h)) &= \notag \\
	-\Big\langle\nabla\cdot\Bigg(\frac{\mathbb{P}_{\mathbb{V}_2}(h_hb_h)\boldsymbol{F}_h}{h_h}\Bigg),\frac{B_h}{h_h}\Big\rangle
	+\Big\langle\frac{B_h}{h_h},\nabla\cdot(b_h\boldsymbol{F}_h)\Big\rangle
	&,
\end{align}
where $\mathbb{P}_{\mathbb{V}_2}(\cdot)$ defines a projection into $\mathbb{V}_2$.
Since the terms do not exactly cancel in either \eqref{eq::tsw_mixed_S1} or \eqref{eq::tsw_mixed_S2}, we do not anticipate 
exact conservation of entropy for this mixed formulation and the magnitude of their difference represents the error in entropy
conservation for the mixed formulation. However numerical experiments detailed in the proceeding section show that loss of entropy 
conservation is small, and solutions remain stable in the absence of dissipation for long times for fully turbulent regimes.
As discussed in Section \ref{sec::tsw_cont}, we note the use of $B_h$, via its application in $\boldsymbol{G}_h$ \eqref{eq::G_from_B}
in the prognostic equation for $b_h$ \eqref{eq::tsw_mixed_disc_b}, and conversely we use $b_h$ in the prognostic equation for $B_h$ 
\eqref{eq::tsw_mixed_disc_B}. This coupling minimises temporal errors associated with the loss of entropy conservation.

As for other compatible discretisations of both the shallow water \cite{MC14,LPG18} and thermal shallow water \cite{Eldred19}
equations, discrete conservation of mass $\mathcal{M}_h^{tsw}=\int h_h\mathrm{d}\Omega^2$, and vorticity, 
$\mathcal{W}_h^{tsw}=\int (h_hq_h-f_h)\mathrm{d}\Omega^2$ is derived by setting $\phi_h = 1$ in \eqref{eq::tsw_mixed_disc_h}
and $\boldsymbol{v}_h=\nabla^{\perp}\psi_h$ in \eqref{eq::tsw_mixed_disc_u} respectively. While both mass and vorticity are 
conserved pointwise for these compatible discretisations of the shallow water equations, for the thermal shallow water
equations vorticity is conserved only globally, owing to the presence of the variable buoyancy in the momentum equation. 
This also holds for the discrete coupled formulation presented below.

\subsubsection{Thermal shallow water equations: coupled formulation}\label{sec::tsw_coupled}

In order to avoid solving separately for both $b_h$ and $B_h$, 
and also recover exact entropy conservation in a discrete form,
we may alternatively construct a discrete analogue of \eqref{eq::tsw_coupled} as
\begin{subequations}\label{eq::tsw_coupled_disc}
	\begin{align}
		\boldsymbol{\mathsf{M}}_1\frac{\mathrm{d}\boldsymbol{u}_h}{\mathrm{d} t} &= -\boldsymbol{\mathsf{C}}_1(q_h,\boldsymbol{F}_h) + 
		\boldsymbol{\mathsf{D}}_2^{\top}\Phi_h^f 
		+ \frac{1}{4}\boldsymbol{\mathsf{M}}_{1*}(b_h',\boldsymbol{\mathsf{M}}_1^{-1}\boldsymbol{\mathsf{D}}_2^{\top}h_h) 
		+ \frac{1}{4}\boldsymbol{\mathsf{D}}_2^{\top}\boldsymbol{\mathsf{M}}_2^{-1}\boldsymbol{\mathsf{M}}_{2*}(b_h',h_h)
		-
		\frac{1}{4}\boldsymbol{\mathsf{M}}_{1*}(\boldsymbol{\mathsf{M}}_1^{-1}\boldsymbol{\mathsf{D}}_2^{\top}b_h',h_h) \\
		\boldsymbol{\mathsf{M}}_2\frac{\mathrm{d} h_h}{\mathrm{d} t} &= -\boldsymbol{\mathsf{D}}_2\boldsymbol{F}_h \\
		\boldsymbol{\mathsf{M}}_2\frac{\mathrm{d} B_h}{\mathrm{d} t} &= -
		\frac{1}{2}\boldsymbol{\mathsf{D}}_2\boldsymbol{\mathsf{M}}_1^{-1}\boldsymbol{\mathsf{M}}_{1*}(b_h',\boldsymbol{F}_h) -
		\frac{1}{2}\boldsymbol{\mathsf{M}}_{2*}(b_h',\boldsymbol{\mathsf{M}}_2^{-1}\boldsymbol{\mathsf{D}}_2\boldsymbol{F}_h) +
		\frac{1}{2}\boldsymbol{\mathsf{K}}_2(\boldsymbol{\mathsf{M}}_1^{-1}\boldsymbol{\mathsf{D}}_2^{\top}b_h',\boldsymbol{F}_h).\label{eq::B_coupled_disc}
	\end{align}
\end{subequations}
Recalling the discrete definition of $b_h'$ \eqref{eq::b_prime},
left 
multiplication of \eqref{eq::tsw_coupled_disc} by $[\boldsymbol{F}_h^{\top}, \Phi_h^{f\top}, T_h^{f\top}]$ leads to the cancellation 
of all right hand side terms, ensuring that there is no spurious generation of energy from the spatial discretisation as
\begin{subequations}
	\begin{align}
		-\boldsymbol{F}_h^{\top}\boldsymbol{\mathsf{C}}_1(q_h,\boldsymbol{F}_h) + 
		\boldsymbol{F}_h^{\top}\boldsymbol{\mathsf{D}}_2^{\top}\Phi_h^f &\notag \\
		+ 
		\frac{1}{4}\boldsymbol{F}_h^{\top}\boldsymbol{\mathsf{M}}_{1*}(b_h',\boldsymbol{\mathsf{M}}_1^{-1}\boldsymbol{\mathsf{D}}_2^{\top}h_h) 
		+
		\frac{1}{4}\boldsymbol{F}_h^{\top}\boldsymbol{\mathsf{D}}_2^{\top}\boldsymbol{\mathsf{M}}_2^{-1}\boldsymbol{\mathsf{M}}_{2*}(b_h',h_h) 
		-
		\frac{1}{4}\boldsymbol{F}_h^{\top}\boldsymbol{\mathsf{M}}_{1*}(\boldsymbol{\mathsf{M}}_1^{-1}\boldsymbol{\mathsf{D}}_2^{\top}b_h',h_h)
		-\Phi_h^{f\top}\boldsymbol{\mathsf{D}}_2\boldsymbol{F}_h &\notag \\
		-\frac{1}{2}T_h^{f\top}\boldsymbol{\mathsf{D}}_2\boldsymbol{\mathsf{M}}_1^{-1}\boldsymbol{\mathsf{M}}_{1*}(b_h',\boldsymbol{F}_h) -
		\frac{1}{2}T_h^{f\top}\boldsymbol{\mathsf{M}}_{2*}(b_h',\boldsymbol{\mathsf{M}}_2^{-1}\boldsymbol{\mathsf{D}}_2\boldsymbol{F}_h) +
		\frac{1}{2}T_h^{f\top}\boldsymbol{\mathsf{K}}_2(\boldsymbol{\mathsf{M}}_1^{-1}\boldsymbol{\mathsf{D}}_2^{\top}b_h',\boldsymbol{F}_h) &= \\
		+ \frac{1}{4}\boldsymbol{F}_h^{\top}\boldsymbol{\mathsf{M}}_{1*}(b_h',\boldsymbol{\mathsf{M}}_1^{-1}\boldsymbol{\mathsf{D}}_2^{\top}h_h) 
		+ \frac{1}{4}\boldsymbol{F}_h^{\top}\boldsymbol{\mathsf{D}}_2^{\top}\boldsymbol{\mathsf{M}}_2^{-1}\boldsymbol{\mathsf{M}}_{2*}(b_h',h_h) 
		- \frac{1}{4}\boldsymbol{F}_h^{\top}\boldsymbol{\mathsf{M}}_{1*}(\boldsymbol{\mathsf{M}}_1^{-1}\boldsymbol{\mathsf{D}}_2^{\top}b_h',h_h)
		&\notag \\
		-\frac{1}{4}h_h^{\top}\boldsymbol{\mathsf{D}}\boldsymbol{\mathsf{M}}_1^{-1}\boldsymbol{\mathsf{M}}_{1*}(b_h',\boldsymbol{F}_h) -
		\frac{1}{4}h_h^{\top}\boldsymbol{\mathsf{M}}_{2*}(b_h',\boldsymbol{\mathsf{M}}_2^{-1}\boldsymbol{\mathsf{D}}_2\boldsymbol{F}_h) +
		\frac{1}{4}h_h^{\top}\boldsymbol{\mathsf{K}}_2(\boldsymbol{\mathsf{M}}_1^{-1}\boldsymbol{\mathsf{D}}_2^{\top}b_h',\boldsymbol{F}_h) &= \\
		\frac{1}{4}\langle\nabla\cdot(b_h'\boldsymbol{F}_h),h_h\rangle -
		\frac{1}{4}\langle\nabla\cdot\boldsymbol{F}_h,h_hb_h'\rangle + 
		\frac{1}{4}\langle\nabla\cdot(h_h'\boldsymbol{F}_h),b_h'\rangle &\notag \\
		-\frac{1}{4}\langle h_h,\nabla\cdot(b_h'\boldsymbol{F}_h)\rangle -
		\frac{1}{4}\langle h_hb_h',\nabla\cdot\boldsymbol{F}_h\rangle + 
		\frac{1}{4}\langle\nabla\cdot(h_h\boldsymbol{F}_h),b_h'\rangle &= 0.
	\end{align}
\end{subequations}

Defining the discrete entropy as $\mathcal{S}^{tsw}_{c,h} = \frac{1}{2}\int h_hb_h'^2\mathrm{d}\Omega^2$, 
the discrete variational derivatives are given as
\begin{subequations}
	\begin{align}
	\boldsymbol{\mathsf{M}}_2\frac{\delta\mathcal{S}_{c,h}^{tsw}}{\delta\boldsymbol{u}_h} &= \boldsymbol{0} \\
	\boldsymbol{\mathsf{M}}_2\frac{\delta\mathcal{S}_{c,h}^{tsw}}{\delta h_h} + 
	\boldsymbol{\mathsf{M}}_2\frac{\delta\mathcal{S}_{c,h}^{tsw}}{\delta b_h'}\cdot\frac{\delta b_h'}{\delta h_h} &= 
	\frac{1}{2}\boldsymbol{\mathsf{M}}_{2*}(b_h',b_h') - 
	\boldsymbol{\mathsf{M}}_2\boldsymbol{\mathsf{M}}_{2*}^{-1}(h_h,
	\boldsymbol{\mathsf{M}}_2\boldsymbol{\mathsf{M}}_{2*}(B_h,
	\boldsymbol{\mathsf{M}}_2\boldsymbol{\mathsf{M}}_{2*}^{-1}(h_h,
	\boldsymbol{\mathsf{M}}_2\boldsymbol{\mathsf{M}}_{2*}(h_h,b_h')))) \notag \\
	&= -\frac{1}{2}\boldsymbol{\mathsf{M}}_{2*}(b_h',b_h')
	\label{eq::dSdh} \\
	\boldsymbol{\mathsf{M}}_2\frac{\delta\mathcal{S}_{c,h}^{tsw}}{\delta B_h} &= \boldsymbol{\mathsf{M}}_{2}\boldsymbol{\mathsf{M}}_{2*}^{-1}(h_h,\boldsymbol{\mathsf{M}}_{2*}(h_h,b_h')) = 
	\boldsymbol{\mathsf{M}}_2b_h', \label{eq::dSdB}
\end{align}\label{eq::tsw_coupled_S_derivs}
\end{subequations}
where in \eqref{eq::dSdB} we have used the discrete identity
$(\delta b_h'/\delta B_h)^{\top} = \boldsymbol{\mathsf{M}}_2\boldsymbol{\mathsf{M}}_{2*}^{-1}(h_h,\cdot)$.
Left multiplication of \eqref{eq::tsw_coupled_disc} by 
the vector of variational derivatives of $\mathcal{S}^{tsw}_{c,h}$ with respect to $\boldsymbol{u}_h$, 
$h_h$, $B_h$ in \eqref{eq::tsw_coupled_S_derivs}
leads to the cancellation of all right hand side terms (again, without
invocation of the chain rule), and hence the conservation of entropy as
\begin{subequations}
	\begin{align}\label{eq::tsw_coupled_S_cons}
		\frac{1}{2}(b_h'^2)^{\top}\boldsymbol{\mathsf{D}}_2\boldsymbol{F}_h 
		-\frac{1}{2}b_h'^{\top}\boldsymbol{\mathsf{D}}_2\boldsymbol{\mathsf{M}}_1^{-1}\boldsymbol{\mathsf{M}}_{1*}(b_h',\boldsymbol{F}_h) -
		\frac{1}{2}b_h'^{\top}\boldsymbol{\mathsf{M}}_{2*}(b_h',\boldsymbol{\mathsf{M}}_2^{-1}\boldsymbol{\mathsf{D}}_2\boldsymbol{F}_h) +
		\frac{1}{2}b_h'^{\top}\boldsymbol{\mathsf{K}}(\boldsymbol{\mathsf{M}}_1^{-1}\boldsymbol{\mathsf{D}}_2^{\top}b_h',\boldsymbol{F}_h) &= \\
	\frac{1}{2}\langle b_h'^2,\nabla\cdot\boldsymbol{F}_h\rangle - \frac{1}{2}\langle b_h',\nabla\cdot(b_h'\boldsymbol{F}_h)\rangle -
	\frac{1}{2}\langle b_h'^2,\nabla\cdot\boldsymbol{F}_h\rangle + \frac{1}{2}\langle\nabla\cdot(b_h'\boldsymbol{F}_h),b_h'\rangle &= 0.
	\end{align}
\end{subequations}
If we once again assume exact integration such that 
$\int\phi_hB_h\mathrm{d}\Omega^2 = \int \phi_hh_hb_h'\mathrm{d}\Omega^2$ yields a unique solution $\forall\phi_h\in\mathbb{V}_2$,
then this conserved entropy can be re-expressed in terms of the prognostic variables as
$\frac{1}{2}\int B_h^2/h_h\mathrm{d}\Omega^2$. Note that unlike the mixed formulation, there is no temporal error in the
relation between $b_h'$ and $B_h$.

As for the energy \cite{BC18,Eldred19,Lee21,LP21}, temporal conservation of the entropy is contingent on the 
preservation of the discrete chain rule as
\begin{equation}
	\frac{\mathrm{d}\mathcal{S}^{tsw}_{c,h}}{\mathrm{d}t} = 
	\Bigg(\frac{\delta\mathcal{S}^{tsw}_{c,h}}{\delta\boldsymbol{u}_h}\Bigg)^{\top}\boldsymbol{\mathsf{M}}_1\frac{d\boldsymbol{u}_h}{dt} + 
	\Bigg(\frac{\delta\mathcal{S}^{tsw}_{c,h}}{\delta h_h}\Bigg)^{\top}\boldsymbol{\mathsf{M}}_2\frac{dh_h}{dt} + 
	\Bigg(\frac{\delta\mathcal{S}^{tsw}_{c,h}}{\delta B_h}\Bigg)^{\top}\boldsymbol{\mathsf{M}}_2\frac{dB_h}{dt}.
\end{equation}
This can only be achieved exactly if the discrete variational derivatives
$\delta\mathcal{S}^{tsw}_{c,h}/\delta\boldsymbol{u}_h$, $\delta\mathcal{S}^{tsw}_{c,h}/\delta h_h$, 
$\delta\mathcal{S}^{tsw}_{c,h}/\delta B_h$ are integrated exactly across the time level to the order 
of the temporal derivative $d/dt$. In our current implementation we use explicit time integration, 
such that the variational derivatives are never evaluated at the new time level, so this is not satisfied
and we instead anticipate some temporal conservation error for both the energy and the entropy.

The coupled formulation above also exactly conserves density weighted buoyancy, $B_h$. However
unlike in the flux form, this is only conserved globally and not locally. Setting the test function 
as $\boldsymbol{1}^{\top}$ in \eqref{eq::B_coupled_disc} we have 
\begin{equation}\label{eq::B_cons_coupled}
	\int\frac{\partial B_h}{\partial t}\mathrm{d}\Omega = 
	-\frac{1}{2}\langle 1_h,\nabla\cdot(b_h'\boldsymbol{F}_h)\rangle 
	-\frac{1}{2}\langle\mathbb{P}_{\mathbb{V}_2}(1_hb_h'),\nabla\cdot\boldsymbol{F}_h\rangle 
	+\frac{1}{2}\langle\nabla\cdot\mathbb{P}_{\mathbb{V}_1}(1_h\boldsymbol{F}_h),b_h'\rangle = 0.
\end{equation}
The first term on the right hand side cancels telescopically due to the compatible mapping of polynomials
in the $\mathbb{V}_1$ space into $\mathbb{V}_2$ via the discrete divergence operator, such that this term 
cancels pointwise, as in the case for the flux form transport equation. The second and third terms on the 
right hand side cancel with one another, assuming that the basis functions that span $\mathbb{V}_2$ are a 
partition of unity, since in that case we have for the projectors into $\mathbb{V}_1$ and $\mathbb{V}_2$ 
respectively that $\mathbb{P}_{\mathbb{V}_1}(1_h\boldsymbol{F}_h) = \boldsymbol{F}_h$ and 
$\mathbb{P}_{\mathbb{V}_2}(1_hb_h') = b_h'$ such that 
$\langle b_h',\nabla\cdot\boldsymbol{F}_h\rangle - \langle\nabla\cdot\boldsymbol{F}_h,b_h'\rangle = 0$.
However for the projector into the $\mathbb{V}_1$ space this is a global operation, and as such 
$B_h$ is conserved only globally and not pointwise for the coupled formulation.

\subsection{3D compressible Euler equations}

Our three dimensional compatible discretisation over the domain $\Omega^3$ is a natural extension 
of the two dimensional discretisation presented in the previous section, for which we introduce the
finite dimensional subspaces $\mathbb{W}_0\subset H^1(\Omega^3)$, 
$\mathbb{W}_1\subset H(\mathrm{curl},\Omega^3)$, $\mathbb{W}_2\subset H(\mathrm{div},\Omega^3)$, 
$\mathbb{W}_3\subset L^2(\Omega^3)$, for which we have a corresponding de Rham complex \cite{AFW06,Cotter23}
\begin{equation}
	\mathbb{R}\longrightarrow\mathbb{W}_0\stackrel{\nabla}{\longrightarrow}
	\mathbb{W}_1\stackrel{\nabla\times}{\longrightarrow}
	\mathbb{W}_2\stackrel{\nabla\cdot}{\longrightarrow}
	\mathbb{W}_3\longrightarrow 0.
\end{equation}
This complex discretely preserves the identities $\nabla\times\nabla :=0$ and $\nabla\cdot\nabla\times :=0$.
As for the two dimensional complex in the previous section, we also recover weak form analogues 
of the \emph{div}, \emph{curl} and \emph{grad} operators as adjoints of the strong form $\nabla$, $\nabla\times$
and $\nabla\cdot$ operators \cite{Natale16,LP21}.

In a slight abuse of notation, we will label the basis functions that span these discrete subspaces
as $\alpha_h\in\mathbb{W}_0$, $\boldsymbol{\psi}_h\in\mathbb{W}_1$, $\boldsymbol{v}_h\in\mathbb{W}_2$, 
$\phi_h\in\mathbb{W}_3$. In a further abuse of notation we will introduce the operators integrated over the
domain $\Omega^3$ as
\begin{subequations}
	\begin{alignat}{5}
		\boldsymbol{\mathsf{M}}_1&:\mathbb{W}_1\rightarrow\mathbb{W}_1'&&\quad\mathrm{for\ which}\quad
		\boldsymbol{\mathsf{M}}_1\boldsymbol{\psi}_h&&:=\langle\boldsymbol{\psi}_h,\cdot\rangle,\quad&&\forall\boldsymbol{\psi}_h\in\mathbb{W}_1,\\
		\boldsymbol{\mathsf{M}}_2&:\mathbb{W}_2\rightarrow\mathbb{W}_2'&&\quad\mathrm{for\ which}\quad
		\boldsymbol{\mathsf{M}}_2\boldsymbol{v}_h&&:=\langle\boldsymbol{v}_h,\cdot\rangle,\quad&&\forall\boldsymbol{v}_h\in\mathbb{W}_2,\\
		\boldsymbol{\mathsf{M}}_3&:\mathbb{W}_3\rightarrow\mathbb{W}_3'&&\quad\mathrm{for\ which}\quad
		\boldsymbol{\mathsf{M}}_3\phi_h&&:=\langle\phi_h,\cdot\rangle,\quad&&\forall\phi_h\in\mathbb{W}_3,\\
		\boldsymbol{\mathsf{D}}_3&:\mathbb{W}_2\rightarrow\mathbb{W}_3'&&\quad\mathrm{for\ which}\quad
		\boldsymbol{\mathsf{D}}_3\boldsymbol{v}_h&&:=\langle\nabla\cdot\boldsymbol{v}_h,\cdot\rangle,\quad&&\forall\boldsymbol{v}_h\in\mathbb{W}_2,\\
		\boldsymbol{\mathsf{R}}_2&:\mathbb{W}_1\rightarrow\mathbb{W}_2'&&\quad\mathrm{for\ which}\quad
		\boldsymbol{\mathsf{R}}_2\boldsymbol{\psi}_h&&:=\langle\nabla\times\boldsymbol{\psi}_h,\cdot\rangle,\quad&&\forall\boldsymbol{\psi}_h\in\mathbb{W}_1,\\
		\boldsymbol{\mathsf{C}}_2&:\mathbb{W}_1\otimes\mathbb{W}_2\rightarrow\mathbb{W}_2'&&\quad\mathrm{for\ which}\quad
		\boldsymbol{\mathsf{C}}_2(\boldsymbol{\psi}_h,\boldsymbol{v}_h)&&:=\langle\boldsymbol{\psi}_h\times\boldsymbol{v}_h,\cdot\rangle,\quad&&\forall\boldsymbol{\psi}_h\in\mathbb{W}_1,\boldsymbol{v}_h\in\mathbb{W}_2,\\
		\boldsymbol{\mathsf{M}}_{1*}&:\mathbb{W}_1\otimes\mathbb{W}_3\rightarrow\mathbb{W}_1'&&\quad\mathrm{for\ which}\quad
		\boldsymbol{\mathsf{M}}_{1*}(\boldsymbol{\psi}_h,\phi_h)&&:=\langle\phi_h\boldsymbol{\psi}_h,\cdot\rangle,\quad&&\forall\boldsymbol{\psi}_h\in\mathbb{W}_1,\phi_h\in\mathbb{W}_3,\\
		\boldsymbol{\mathsf{M}}_{2*}&:\mathbb{W}_2\otimes\mathbb{W}_3\rightarrow\mathbb{W}_3'&&\quad\mathrm{for\ which}\quad
		\boldsymbol{\mathsf{M}}_{2*}(\boldsymbol{v}_h,\phi_h)&&:=\langle\phi_h\boldsymbol{v}_h,\cdot\rangle,\quad&&\forall\boldsymbol{v}_h\in\mathbb{W}_2,\phi_h\in\mathbb{W}_3,\\
		\boldsymbol{\mathsf{M}}_{3*}&:\mathbb{W}_3\otimes\mathbb{W}_3\rightarrow\mathbb{W}_3'&&\quad\mathrm{for\ which}\quad
		\boldsymbol{\mathsf{M}}_{3*}(\phi_h,\phi_k)&&:=\langle\phi_h\phi_k,\cdot\rangle,\quad&&\forall\phi_h,\phi_k\in\mathbb{W}_3,\\
		\boldsymbol{\mathsf{K}}_3&:\mathbb{W}_2\otimes\mathbb{W}_2\rightarrow\mathbb{W}_3'&&\quad\mathrm{for\ which}\quad
		\boldsymbol{\mathsf{K}}_3(\boldsymbol{v}_h,\boldsymbol{v}_k)&&:=\langle\boldsymbol{v}_h\cdot\boldsymbol{v}_k,\cdot\rangle,\quad&&\forall\boldsymbol{v}_h,\boldsymbol{v}_k\in\mathbb{W}_2,\\
		\boldsymbol{\mathsf{D}}_{3*}&:\mathbb{W}_2\otimes\mathbb{W}_3\rightarrow\mathbb{W}_3'&&\quad\mathrm{for\ which}\quad
		\boldsymbol{\mathsf{D}}_{3*}(\boldsymbol{v}_h,\phi_h)&&:=\langle\nabla\cdot\boldsymbol{v}_h,\phi_h\cdot\rangle,\quad&&\forall\boldsymbol{v}_h\in\mathbb{W}_2,\phi_h\in\mathbb{W}_3.
	\end{alignat}
\end{subequations}

Having defined the various operators, the discrete potential vorticity, $\boldsymbol{q}_h\in\mathbb{W}_1$ and 
variational derivatives $\boldsymbol{U}_h\in\mathbb{W}_2$, $\Phi_h, \Pi_h\in\mathbb{W}_3$ are given as
\begin{subequations}
	\begin{align}
		\boldsymbol{\mathsf{M}}_{1*}(\boldsymbol{q}_h,h_h) &= -\boldsymbol{\mathsf{R}}_2^{\top}\boldsymbol{u}_h + 
		\langle\boldsymbol{\psi}_h,f_h\boldsymbol{e}_z\rangle \\
		\boldsymbol{\mathsf{M}}_2\boldsymbol{U}_h &= \boldsymbol{\mathsf{M}}_{2*}(\boldsymbol{u}_h,\rho_h) \\
		\boldsymbol{\mathsf{M}}_3\Phi_h &= \frac{1}{2}\boldsymbol{\mathsf{K}}_3(\boldsymbol{u}_h,\boldsymbol{u}_h) + 
		\langle\phi_h,gz\boldsymbol{e}_z\rangle \\
		\boldsymbol{\mathsf{M}}_3\Pi_h &= c_p\Big(\frac{R}{p_0}\Big)^{\frac{R}{c_v}}\Big\langle\phi_h,\Theta_h^{\frac{R}{c_v}}\Big\rangle,
	\end{align}
\end{subequations}
where $\boldsymbol{e}_z$ is the unit vector in the vertical direction. Similarly to \eqref{eq::b_prime}, we define the 
discrete potential temperature, $\theta_h\in\mathbb{W}_3$ via
\begin{equation}\label{eq::Theta_diag}
	\boldsymbol{\mathsf{M}}_{3*}(\rho_h,\theta_h) = \boldsymbol{\mathsf{M}}_3\Theta_h.
\end{equation}
The semi-discrete equivalent of the prognostic equations \eqref{eq::ce_coupled} for 
$\boldsymbol{u}_h\in\mathbb{W}_2$, $\rho_h, \Theta_h\in\mathbb{W}_3$ are then expressed as
\begin{subequations}\label{eq::ce_disc}
	\begin{align}
		\boldsymbol{\mathsf{M}}_2\frac{\mathrm{d}\boldsymbol{u}_h}{\mathrm{d} t} &= -\boldsymbol{\mathsf{C}}_2(\boldsymbol{q}_h,\boldsymbol{U}_h) + 
		\boldsymbol{\mathsf{D}}_3^{\top}\Phi_h 
		+ \frac{1}{2}\boldsymbol{\mathsf{M}}_{2*}(\theta_h,\boldsymbol{\mathsf{M}}_2^{-1}\boldsymbol{\mathsf{D}}_3^{\top}\Pi_h) 
		+ \frac{1}{2}\boldsymbol{\mathsf{D}}_{3*}^{\top}(\theta_h,\Pi_h) 
		- \frac{1}{2}\boldsymbol{\mathsf{M}}_{2*}(\boldsymbol{\mathsf{M}}_2^{-1}\boldsymbol{\mathsf{D}}_3^{\top}\theta_h,\Pi_h) \\
		\boldsymbol{\mathsf{M}}_3\frac{\mathrm{d}\rho_h}{\mathrm{d} t} &= -\boldsymbol{\mathsf{D}}_3\boldsymbol{U}_h \\
		\boldsymbol{\mathsf{M}}_3\frac{\mathrm{d}\Theta_h}{\mathrm{d} t} &= -
		\frac{1}{2}\boldsymbol{\mathsf{D}}_3\boldsymbol{\mathsf{M}}_2^{-1}\boldsymbol{\mathsf{M}}_{2*}(\theta_h,\boldsymbol{U}_h) -
		\frac{1}{2}\boldsymbol{\mathsf{M}}_{3*}(\theta_h,\boldsymbol{\mathsf{M}}_3^{-1}\boldsymbol{\mathsf{D}}_3\boldsymbol{U}_h) +
		\frac{1}{2}\boldsymbol{\mathsf{K}}_3(\boldsymbol{\mathsf{M}}_2^{-1}\boldsymbol{\mathsf{D}}_3^{\top}\theta_h,\boldsymbol{U}_h).
	\end{align}
\end{subequations}

As for the coupled formulation of the thermal shallow water equations, energy conservation is
assured via the left-multiplication of \eqref{eq::ce_disc} by $[\boldsymbol{U}_h^{\top}, \Phi_h^{\top}, \Pi_h^{\top}]$, such that
the right hand side cancels as
\begin{subequations}
	\begin{align}
		-\boldsymbol{U}_h^{\top}\boldsymbol{\mathsf{C}}_2(\boldsymbol{q}_h,\boldsymbol{U}_h) + 
		\boldsymbol{U}_h^{\top}\boldsymbol{\mathsf{D}}_3^{\top}\Phi_h 
		+ \frac{1}{2}\boldsymbol{U}_h^{\top}\boldsymbol{\mathsf{M}}_{2*}(\theta_h,\boldsymbol{\mathsf{M}}_2^{-1}\boldsymbol{\mathsf{D}}_3^{\top}\Pi_h) 
		+ \frac{1}{2}\boldsymbol{U}_h^{\top}\boldsymbol{\mathsf{D}}_{3*}^{\top}(\theta_h,\Pi_h) 
		-
		\frac{1}{2}\boldsymbol{U}_h^{\top}\boldsymbol{\mathsf{M}}_{2*}(\boldsymbol{\mathsf{M}}_2^{-1}\boldsymbol{\mathsf{D}}_3^{\top}\theta_h,\Pi_h) &\notag \\
		-\Phi_h^{\top}\boldsymbol{\mathsf{D}}_3\boldsymbol{U}_h
		-\frac{1}{2}\Pi_h^{\top}\boldsymbol{\mathsf{D}}_3\boldsymbol{\mathsf{M}}_2^{-1}\boldsymbol{\mathsf{M}}_{2*}(\theta_h,\boldsymbol{U}_h) -
		\frac{1}{2}\Pi_h^{\top}\boldsymbol{\mathsf{M}}_{3*}(\theta_h,\boldsymbol{\mathsf{M}}_3^{-1}\boldsymbol{\mathsf{D}}_3\boldsymbol{U}_h) +
		\frac{1}{2}\Pi_h^{\top}\boldsymbol{\mathsf{K}}_3(\boldsymbol{\mathsf{M}}_2^{-1}\boldsymbol{\mathsf{D}}_3^{\top}\theta_h,\boldsymbol{U}_h) &= \\
		\frac{1}{2}\langle\nabla\cdot(\theta_h\boldsymbol{U}_h),\Pi_h\rangle
		+ 
		\frac{1}{2}\langle\nabla\cdot\boldsymbol{U}_h,\theta_h\Pi_h\rangle 
		-
		\frac{1}{2}\langle\nabla\cdot(\Pi_h\boldsymbol{U}_h),\theta_h\rangle &\notag \\
		-\frac{1}{2}\langle\Pi_h,\nabla\cdot(\theta_h\boldsymbol{U}_h)\rangle -
		\frac{1}{2}\langle\theta_h\Pi_h,\nabla\cdot\boldsymbol{U}_h\rangle + \frac{1}{2}\langle\nabla\cdot(\Pi_h\boldsymbol{U}_h),\theta_h\rangle &= 0.
	\end{align}
\end{subequations}
Conservation of entropy is established in an almost identical fashion as for the discrete coupled thermal shallow
water equations in \eqref{eq::tsw_coupled_S_cons}. For the entropy $\mathcal{S}^{ce}_{h} = \frac{1}{2}\int\rho_h\theta_h^2\mathrm{d}\Omega^3$, 
left multiplication of \eqref{eq::ce_disc} by 
$(\nabla\mathcal{S}^{ce}_h)^{\top} = [\boldsymbol{0}^{\top}, -\frac{1}{2}\theta_h^{2\top}, \theta_h^{\top}]$,
which are determined in a direct analogue of \eqref{eq::tsw_coupled_S_derivs}, gives
\begin{subequations}
	\begin{align}
		\frac{1}{2}\theta_h^{2\top}\boldsymbol{\mathsf{D}}_3\boldsymbol{U}_h 
		-\frac{1}{2}\theta_h^{\top}\boldsymbol{\mathsf{D}}_3\boldsymbol{\mathsf{M}}_2^{-1}\boldsymbol{\mathsf{M}}_{2*}(\theta_h,\boldsymbol{U}_h) -
		\frac{1}{2}\theta_h^{\top}\boldsymbol{\mathsf{M}}_{3*}(\theta_h,\boldsymbol{\mathsf{M}}_3^{-1}\boldsymbol{\mathsf{D}}_3\boldsymbol{U}_h) +
		\frac{1}{2}\theta_h^{\top}\boldsymbol{\mathsf{K}}_3(\boldsymbol{\mathsf{M}}_2^{-1}\boldsymbol{\mathsf{D}}_3^{\top}\theta_h,\boldsymbol{U}_h) &= \\
		\frac{1}{2}\langle\theta_h^2,\nabla\cdot\boldsymbol{U}_h\rangle -
		\frac{1}{2}\langle\theta_h,\nabla\cdot(\theta_h\boldsymbol{U}_h)\rangle -
		\frac{1}{2}\langle\theta_h^2,\nabla\cdot\boldsymbol{U}_h\rangle +
		\frac{1}{2}\langle\theta_h,\nabla\cdot(\theta_h\boldsymbol{U}_h)\rangle &= 0.
	\end{align}
\end{subequations}

\section{Results}

\subsection{Thermal shallow water equations: thermogeostrophic balance}

In order to first demonstrate the correctness of the entropy conserving construction, the thermal 
shallow water equations are solved in both mixed \eqref{eq::tsw_mixed_disc} and coupled 
\eqref{eq::tsw_coupled_disc} form using a cubed sphere discretisation \cite{LP18}
for an analytic test case for thermogeostropic balance. This represents an extension
of a similar planar test case \cite{Eldred19} extended to the sphere \cite{Ricardo23b} for a steady
solution of the form 
\begin{subequations}
	\begin{align}
		u(\lambda,\chi) &= u_0\cos(\chi) \\
		v(\lambda,\chi) &= 0 \\
		h(\lambda,\chi) &= h_0 - \frac{u_0}{g}\Big(r_e\omega_e + \frac{1}{2}u_0\Big)\sin^2(\chi) \\
		b(\lambda,\chi) &= g\Bigg(1 + 0.05\Bigg(\frac{h_0}{h}\Bigg)^2\Bigg),
	\end{align}
\end{subequations}
where $u_0 = 2\pi r_e/T$, with $r_e = 6371220.0m$ being the earth's radius and $T=12\ days$, 
$h_0 = 2.94\times 10^4/g$, $\omega_e=7.292\times 10^{-5}s^{-1}$ being the earth's angular 
frequency and $(\lambda,\chi)$ being the longitudinal and latitudinal coordinates respectively
in the physical domain $\Omega^2$.

As for all test cases presented here, we use third order elements and six point Gauss-Lobatto-Legendre
quadrature, which is exact for polynomials of degree nine and lower. However since our Jacobian 
involves trigonometric functions \cite{Guba14,LP18} which cannot be integrated exactly, we do not 
preserve the properties of exact spatial integration in our results.

\begin{figure}[!hbtp]
\begin{center}
\includegraphics[width=0.32\textwidth,height=0.24\textwidth]{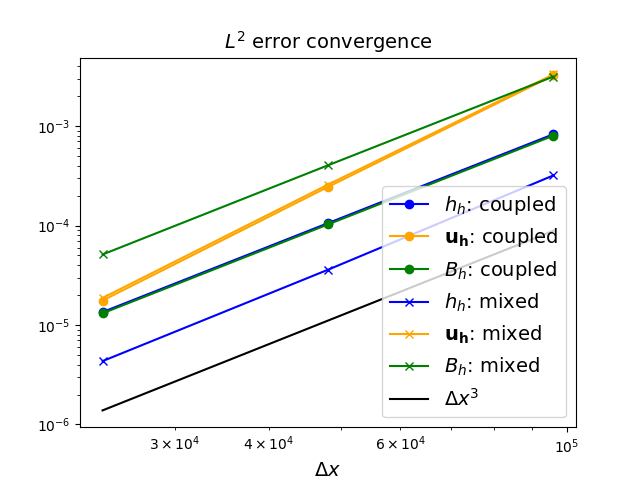}
\includegraphics[width=0.32\textwidth,height=0.24\textwidth]{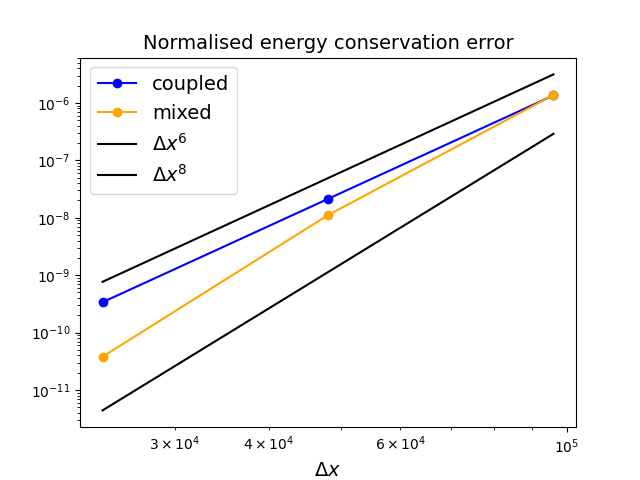}
\includegraphics[width=0.32\textwidth,height=0.24\textwidth]{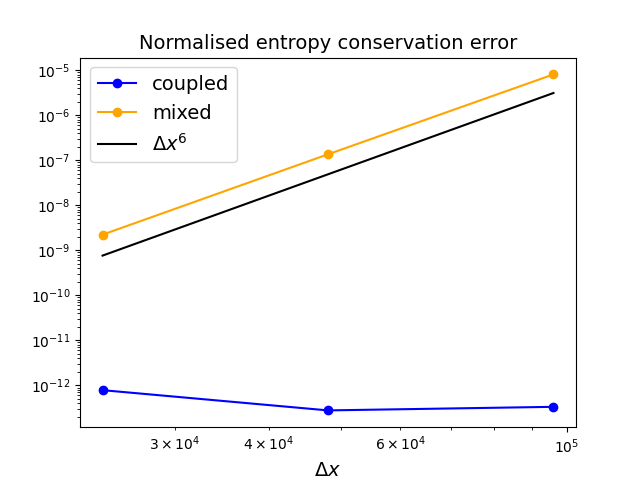}
\caption{$L^2$ errors (left), energy conservation error (center) and entropy conservation error (right) 
convergence at day 5 of the thermogeostrophic test case using the mixed and coupled formulations 
of the thermal shallow water equations.}
\label{fig::tsw_convergence}
\end{center}
\end{figure}

Figure \ref{fig::tsw_convergence} shows the convergence of errors in the $L^2$ norm for $h_h$,
$\boldsymbol{u}_h$ and $B_h$ for both the mixed \eqref{eq::tsw_mixed_disc} and coupled \eqref{eq::tsw_coupled_disc}
formulations (in both cases using a stiffly-stable third order Runge-Kutta time integrator \cite{SO88})
after 5 days of integration, using $6\times 8^2$, $6\times 16^2$ and $6\times 32^2$ third order
elements on the cubed sphere (and corresponding time steps of $\Delta t=120s$, $\Delta t=60s$
and $\Delta t=30s$). In all cases the errors converge with the anticipated third power
of the average grid spacing, $\Delta x$, with the exception of the velocities, $\boldsymbol{u}_h$,
which exhibit a steeper rate of convergence due to the higher polynomial degree in the normal 
direction for the basis functions spanning the $\mathbb{V}_1\subset H(\mathrm{div},\Omega^2)$
subspace \cite{LP18}. 

This figure also shows the energy and entropy conservation errors for both formulations. In terms 
of energy, this converges at a higher rate for the mixed formulation, with the temporal errors being 
associated with the inexact representation of the temporal chain rule \eqref{eq::temporal_en_cons}.
As mentioned previously, energy can be conserved in time also if the temporal chain rule is preserved
exactly using a semi-implicit time integrator with the variational derivatives integrated exactly
in time \cite{BC18,Eldred19,Lee21,LP21}. For the entropy, this is conserved at close to machine 
precision for the coupled scheme, whereas the mixed scheme involves a (convergent) temporal error
due to the independent time stepping of both $b_h$ and $B_h$, and the loss of exact equality in the
cancellation of the discrete terms in \eqref{eq::tsw_mixed_S1} and \eqref{eq::tsw_mixed_S2}.

\subsection{Thermal shallow water equations: shear flow instability}

In order to study the conservation properties of the two different energy and entropy conserving 
thermal shallow water formulations in a more dynamically rich setting we present results for a 
modified shear flow test instability case \cite{Galewsky04}. 
In the original test case for the rotating shallow water equations with constant gravity a 
geostrophically balanced zonal jet is overlaid with a Gaussian perturbation in the depth field 
which generates a gravity wave that triggers an instability in the jet. 
For the thermal shallow water equations with variable buoyancy we extended this test case
with the addition of a perturbed buoyancy field of the form 
\begin{equation}
	b(\lambda,\chi) = g\Big(1 - 0.1\cos(\chi)e^{-(3\lambda)^2-(15(\pi/4-\chi))^2}\Big),
\end{equation}
which is similar to the Gaussian perturbation to the otherwise geostrophically balanced depth field
in the original shallow water configuration, which we also include in our configuration. For both
the mixed and coupled formulations the model was configured with $6\times 32^2$ third order elements and a time step of
$\Delta t=30s$, and run for 20 days so as to ensure model stability for a mature turbulent state.
No energy or entropy damping is applied in either case. Also we have also avoided the use of any 
potential enstrophy damping, as is typically required to stabilise energy conserving discretisations
of the constant buoyancy shallow water equations \cite{SB85,Lee22}.

\begin{figure}[!hbtp]
\begin{center}
\includegraphics[width=0.32\textwidth,height=0.24\textwidth]{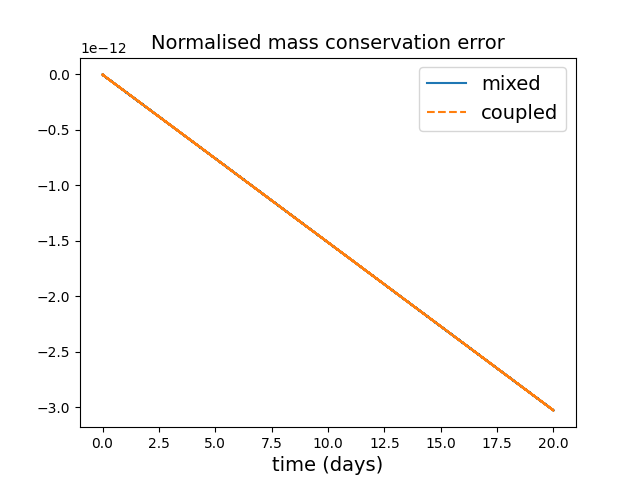}
\includegraphics[width=0.32\textwidth,height=0.24\textwidth]{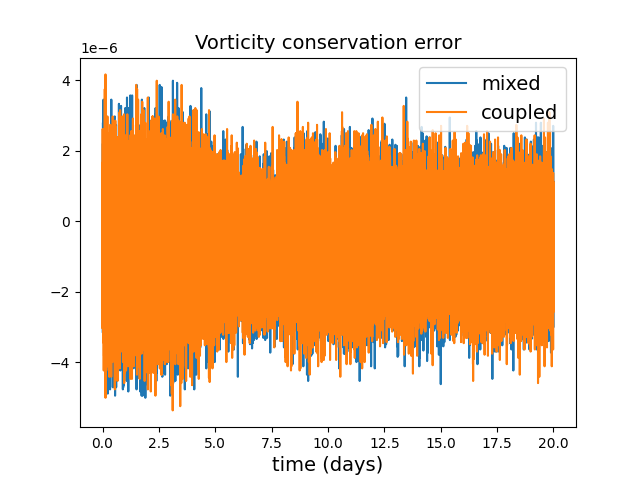}
\includegraphics[width=0.32\textwidth,height=0.24\textwidth]{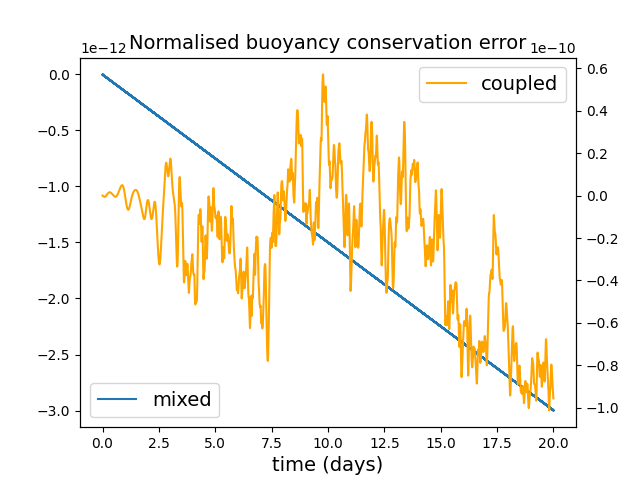}
\caption{Conservation errors for the mixed and coupled formulations of the thermal shallow water 
	equations for the mass (left), un-normalised vorticity (center) and buoyancy (right) for the shear flow instability test case.}
\label{fig::tsw_conservation_1}
\end{center}
\end{figure}

Figure \ref{fig::tsw_conservation_1} shows the conservation errors
for the shear flow test case for the mass, vorticity and depth weighted buoyancy. Both the mixed
and coupled formulations exhibit a machine precision loss of mass at each time step, as has been
observed previously for third order Runge-Kutta time integration \cite{LP20}.
The un-normalised vorticity conservation errors are consistent with those observed for previous compatible finite
element discretisations of the rotating shallow water equations \cite{LP18,Lee22} for the shear flow 
instability test case \cite{Galewsky04}. While the depth weighted buoyancy conservation errors
are of comparable magnitude to the mass conservation errors for the mixed formulation, these
are $\mathcal{O}(100)$ times larger for the coupled formulation for which density weighted
buoyancy is conserved only globally and not locally as described in Section \ref{sec::tsw_coupled}, 
and for which conservation is dependent on the accuracy of the linear solve that represents the 
projection into the $\mathbb{V}_1$ space as given in \eqref{eq::B_cons_coupled}.

\begin{figure}[!hbtp]
\begin{center}
\includegraphics[width=0.48\textwidth,height=0.36\textwidth]{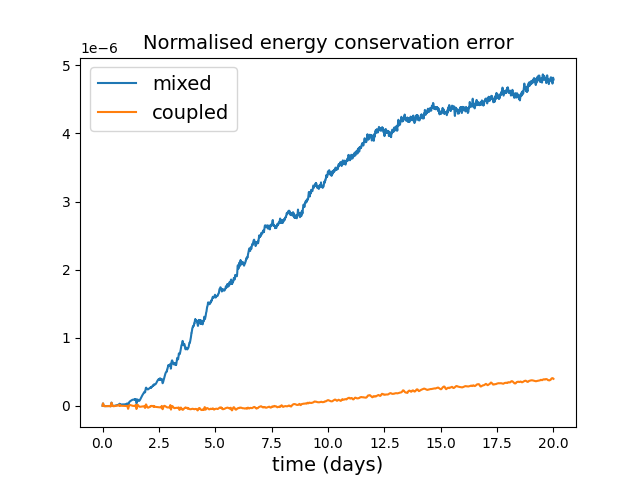}
\includegraphics[width=0.48\textwidth,height=0.36\textwidth]{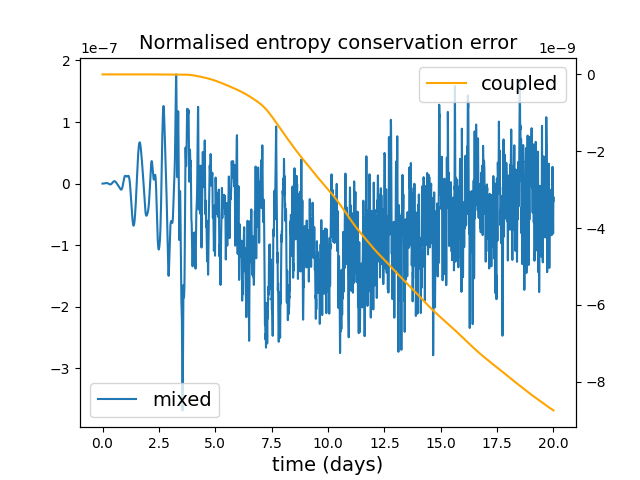}
\caption{Conservation errors for the mixed and coupled formulations of the thermal shallow water 
	equations for the energy (left) and entropy (right) for the shear flow instability test case.}
\label{fig::tsw_conservation_2}
\end{center}
\end{figure}

The energy and entropy conservation errors for the two formulations are presented in Fig.
\ref{fig::tsw_conservation_2}. 
While the energy conservation errors are $\mathcal{O}(10)$ times larger for the mixed formulation,
the entropy conservation errors are $\mathcal{O}(100)$ times larger, perhaps due to numerical 
errors in the cancellation of terms in \eqref{eq::tsw_mixed_S1} and \eqref{eq::tsw_mixed_S2}.
As discussed in Section \ref{sec::tsw_coupled}, since we are using explicit time integration for
which the variational derivatives of the energy and entropy are not exactly integrated in time
to the order of the time stepping scheme, we anticipate some temporal loss of conservation in
both these quantities.
Nevertheless with the application of the strong stability preserving explicit integrator
both formulations are stable for the duration of the simulation without any form of stabilisation.

\begin{figure}[!hbtp]
\begin{center}
\includegraphics[width=0.48\textwidth,height=0.36\textwidth]{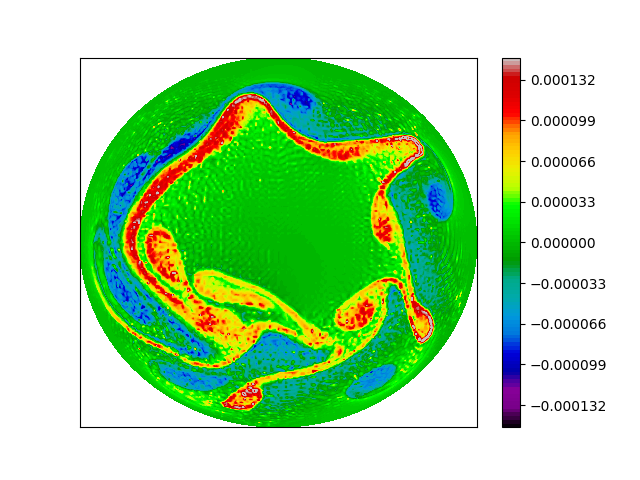}
\includegraphics[width=0.48\textwidth,height=0.36\textwidth]{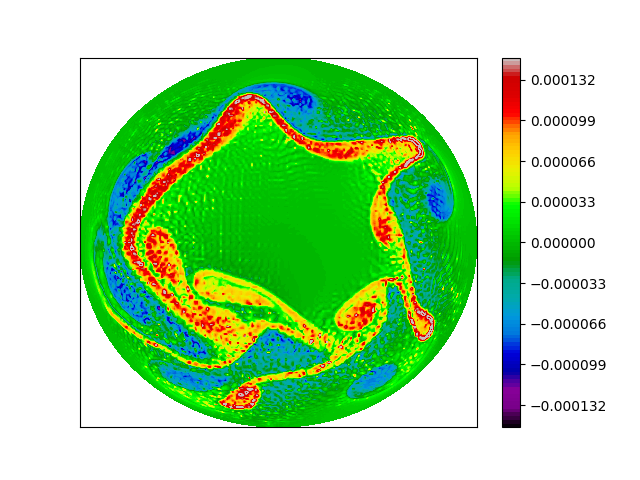}
\caption{Vorticity field (northern hemisphere) at day 6 for the thermal shallow water 
	equations using the mixed (left) and coupled (right) formulations for the shear flow instability test case.
	\blue{Values are given in units of $s^{-1}$}}
\label{fig::tsw_vorticity}
\end{center}
\end{figure}

\begin{figure}[!hbtp]
\begin{center}
\includegraphics[width=0.48\textwidth,height=0.36\textwidth]{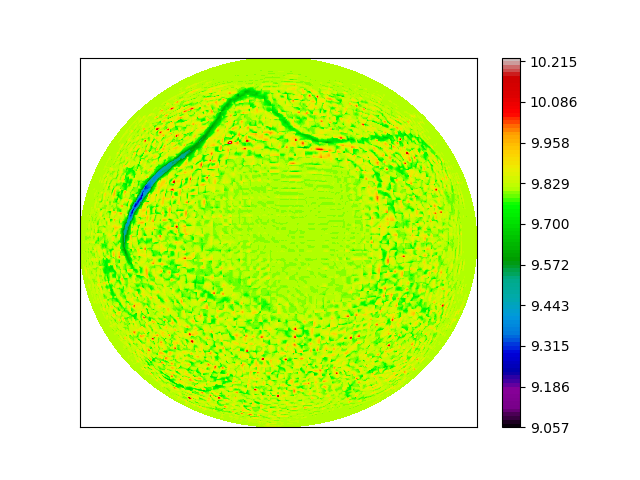}
\includegraphics[width=0.48\textwidth,height=0.36\textwidth]{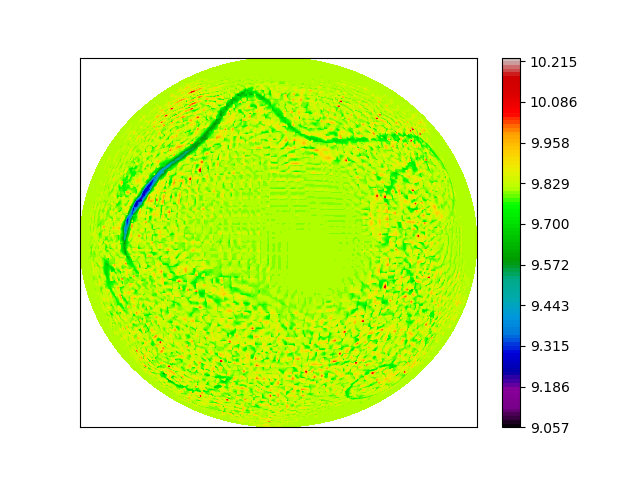}
\caption{Buoyancy field (northern hemisphere) at day 6 for the thermal shallow water 
	equations using the mixed (left) and coupled (right) formulations for the shear flow instability test case.
	\blue{Values are given in units of $m\cdot s^{-2}$.}}
\label{fig::tsw_buoyancy}
\end{center}
\end{figure}

The northern-hemisphere vorticity and buoyancy fields are shown at day 6 for both schemes
in Figs. \ref{fig::tsw_vorticity} and \ref{fig::tsw_buoyancy}. The results are similar for
both schemes. Since no energy or entropy damping is applied, both develop a large amount
of grid scale noise due to nonlinear aliasing. However since both energy and entropy are conserved by the 
spatial discretisation, this noise does not trigger any instability in the model, and both 
are able to run stably for the duration of the simulation.

\subsection{3D compressible Euler equations: baroclinic instability}
We validate the entropy and energy conserving formulation for the 3D compressible
Euler equations \eqref{eq::ce_disc} for a standard baroclinic instability test case 
on the sphere \cite{UMJS14}, using $6\times 24^2$ elements of third order in the 
horizontal dimensions and first order in the vertical dimension, with 30 vertical 
levels and a time step of $\Delta t=75s$ for 12 days. The model uses a second order 
horizontally explicit/vertically implicit time splitting that exactly balances all
energy exchanges \cite{LP21}, and an exact energy conserving implicit solver for 
the vertical \cite{Lee21}. 
While the preconditioner for this implicit solver was originally designed 
for the flux form of the density weighted potential temperature, and the corresponding
energy conserving pressure gradient term, it also provides robust convergence for the new 
coupled energy entropy conserving formulation.
The model uses a biharmonic viscosity in the horizontal for the momentum equation \cite{Guba14},
\blue{which is observed to be necessary for model stability. Both the shallow water \cite{Lee22},
and thermal shallow water equations (present article) are observed to be stable in the absence of
dissipation or upwinding if the higher order invariants are conserved (energy and potential enstrophy for the shallow
water equations, and energy and entropy for the thermal shallow water equations). Therefore we reason that 
in the present case biharmonic viscosity is required for the 3D compressible Euler equations 
due to splitting errors in the horizontally explicit, vertically implicit time integrator \cite{LP21}, 
and the lack of temporal conservation of either energy (which in three dimensions will cascade to small scales) 
or entropy in the explicit horizontal integrator}. The code also uses a Rayleigh 
damping term on the top three levels in the vertical momentum equation in order to 
suppress vertical oscillations associated with the hydrostatic adjustment at the 
beginning of the simulation. 

The model is compared to a near identical formulation \cite{Lee21,LP21}, for which 
the density weighted potential temperature is transported in flux form 
only, with a corresponding energy conserving adjoint form of the pressure gradient 
in the momentum equation. This original flux formulation requires the inclusion of an additional biharmonic viscosity term 
for the flux form density weighted potential temperature equation in order to ensure 
long term stability. In the present case neither the flux form nor the coupled form 
is run with biharmonic viscosity or any other form of stabilisation on the density 
weighted potential temperature equation so as to demonstrate the loss of entropy stability 
for the original flux formulation. While the new coupled entropy conserving formulation is
able to run stably for the duration of the simulation with no sign of instability, the
original flux formulation goes unstable due to unbounded entropy conservation errors
after 11 days of simulation time.

\begin{figure}[!hbtp]
\begin{center}
\includegraphics[width=0.32\textwidth,height=0.24\textwidth]{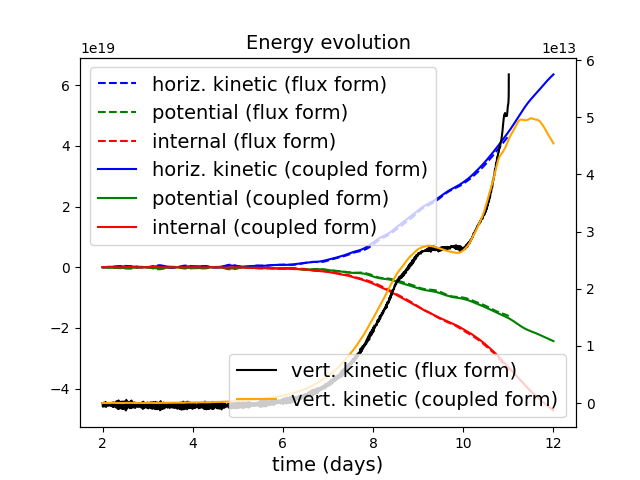}
\includegraphics[width=0.32\textwidth,height=0.24\textwidth]{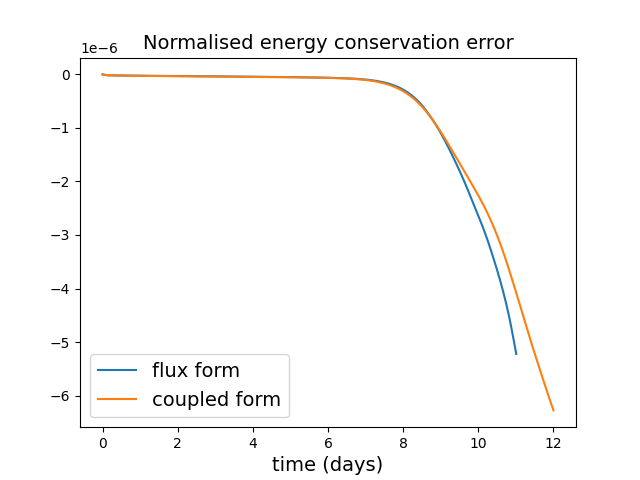}
\includegraphics[width=0.32\textwidth,height=0.24\textwidth]{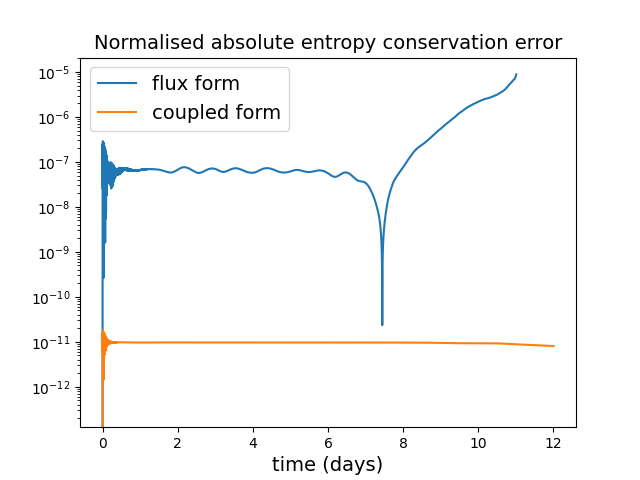}
\caption{Temporal evolution of the energy components for the 3D compressible Euler
equations for the flux and coupled form potential temperature formulations (left),
	and the total energy (center) and absolute value entropy (right) 
	conservation errors for the baroclinic instability test case.}
\label{fig::ce_energetics_1}
\end{center}
\end{figure}

Figure \ref{fig::ce_energetics_1} shows the evolution of the kinetic (vertical and 
horizontal), potential and internal energies (difference from initial values) for the
flux form and coupled formulations. Note the difference in scale for the vertical 
kinetic energy evolution (right hand side axis on the first figure), for which the 
energy changes are $\mathcal{O}(10^6)$ times smaller than for the other energy components.
While the horizontal kinetic, potential and internal energy differences are negligible
between the two formulations (with slightly less observable oscillation for the coupled scheme), 
there is a sharp uptick in the vertical kinetic energy
for the flux form shortly before the simulation goes unstable. This is not the case
for the coupled entropy conserving formulation, which remains stable while accurately
capturing the evolution of the vertical kinetic energy due to the baroclinic process.

The instability of the flux formulation is preceeded by a steady growth in entropy,
which is not observed for the entropy conserving coupled formulation, for which the 
entropy conservation error remains $\mathcal{O}(10^4)$ times smaller than the flux 
form for the duration of the simulation (but is still not machine precision due to the loss
of entropy conservation in the temporal discretisation). This is despite
the fact that the coupled formulation also exhibits a slower decay in energy, presumably
since there is less grid scale noise to be smoothed out by the biharmonic viscosity
on the momentum equation.

\begin{figure}[!hbtp]
\begin{center}
\includegraphics[width=0.32\textwidth,height=0.24\textwidth]{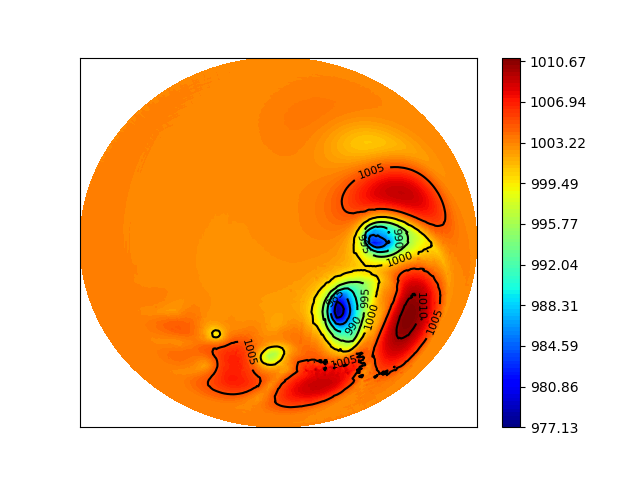}
\includegraphics[width=0.32\textwidth,height=0.24\textwidth]{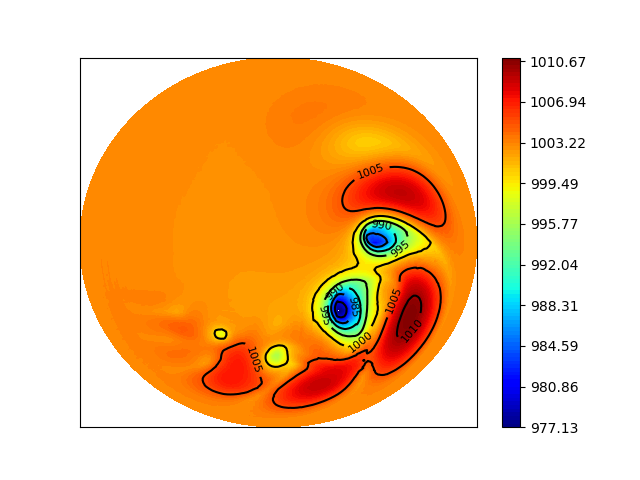}
\includegraphics[width=0.32\textwidth,height=0.24\textwidth]{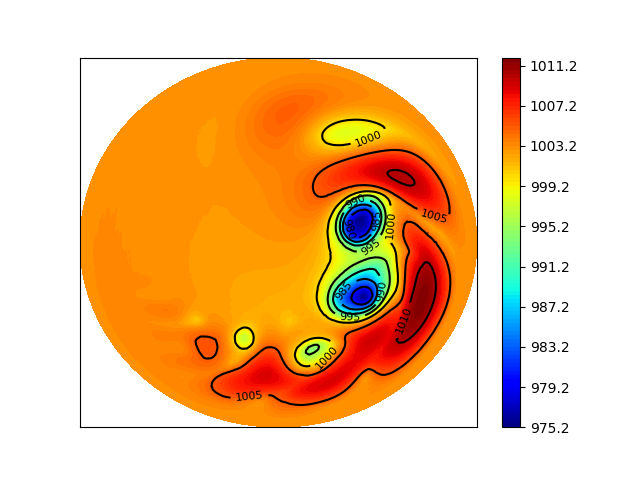}
	\caption{Surface level Exner pressure for the original (flux) formulation and entropy conserving
	(coupled) formulation at day 11 (left and center), and entropy conserving formulation at day 12 
	(right) for the baroclinic instability test case. \blue{Values are given in units of Joules per 
	kilogram per Kelvin, $J\cdot kg^{-1}\cdot K^{-1}$.}}
\label{fig::ce_exner}
\end{center}
\end{figure}

The lowest level Exner pressure is presented at day 11 for both the flux form and coupled form
simulations, as well as for the coupled form at day 12 (after the flux form has gone unstable) 
in Fig. \ref{fig::ce_exner}. \blue{Since the Exner pressure is scaled by the specific heat at
constant pressure in \eqref{eq::eos}, this figure is presented in units of $J\cdot kg^{-1}\cdot K^{-1}$}. 
While the onset of numerical instability for the flux form is visible 
from the somewhat noisy structure of the Exner pressure contours, this is not the case for the 
coupled form, which remain smooth and coherent after 12 days of simulation.

\begin{figure}[!hbtp]
\begin{center}
\includegraphics[width=0.32\textwidth,height=0.24\textwidth]{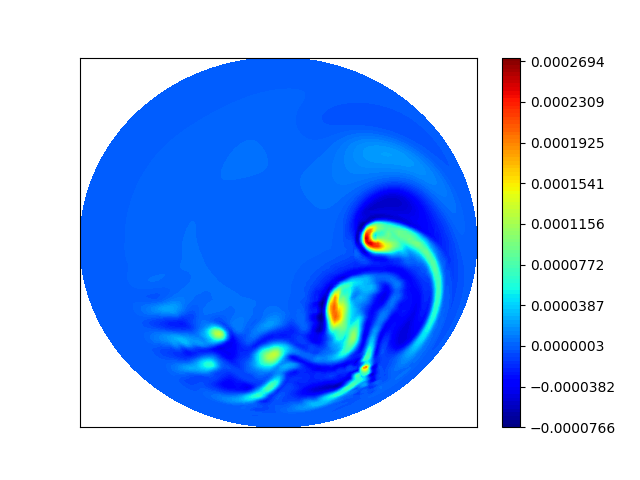}
\includegraphics[width=0.32\textwidth,height=0.24\textwidth]{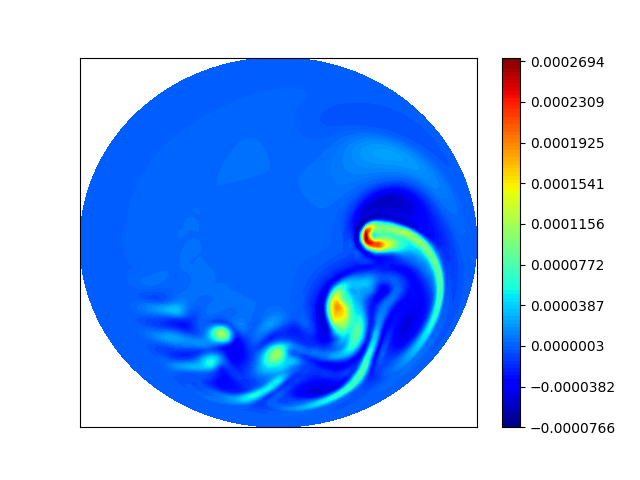}
\includegraphics[width=0.32\textwidth,height=0.24\textwidth]{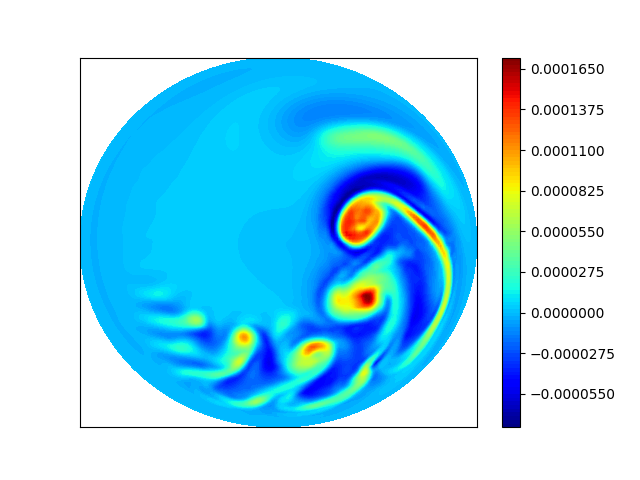}
	\caption{Vertical vorticity component at $z\approx 1.5km$ for the original (flux) formulation 
	and entropy conserving (coupled) formulation at day 11 (left and center), and entropy conserving 
	formulation at day 12 (right) for the baroclinic instability test case. \blue{Values are
	given in units of inverse time, $s^{-1}$.}}
\label{fig::ce_vorticity}
\end{center}
\end{figure}

\begin{figure}[!hbtp]
\begin{center}
\includegraphics[width=0.32\textwidth,height=0.24\textwidth]{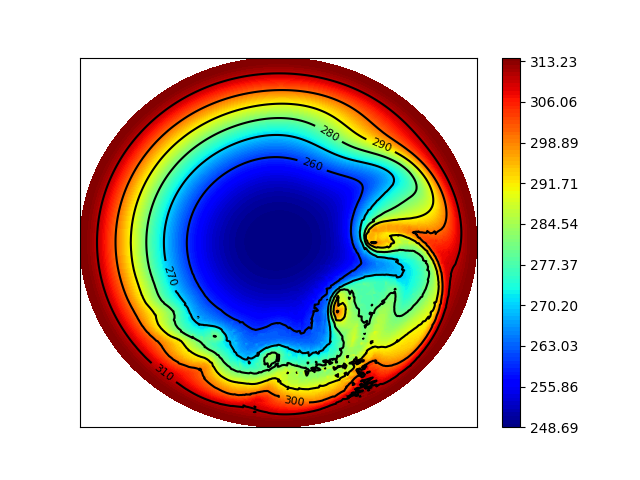}
\includegraphics[width=0.32\textwidth,height=0.24\textwidth]{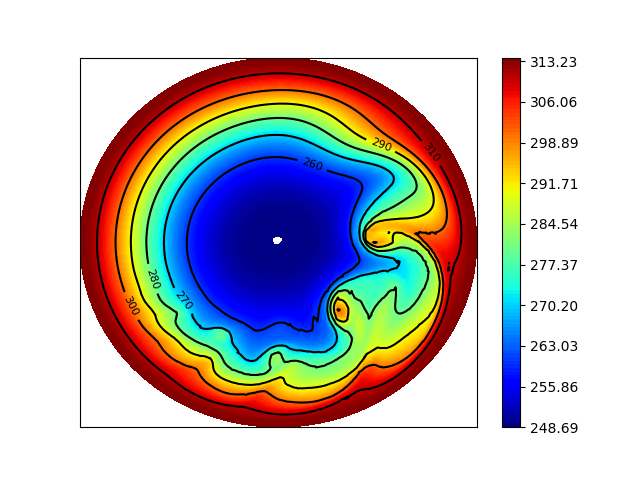}
\includegraphics[width=0.32\textwidth,height=0.24\textwidth]{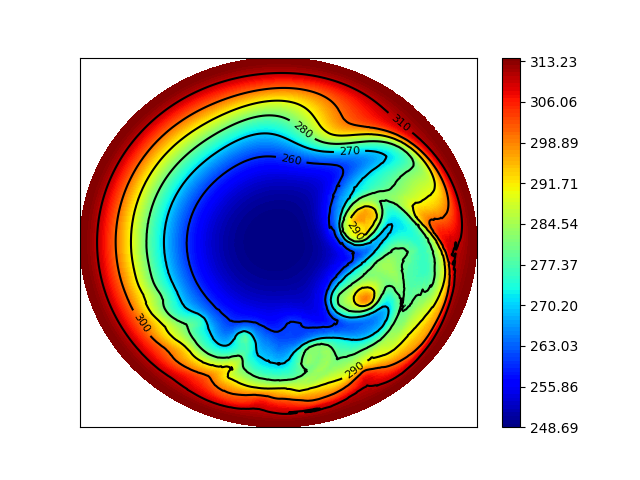}
	\caption{Potential temperature at $z\approx 1.25km$ for the original (flux) formulation 
	and entropy conserving (coupled) formulation at day 11 (left and center), and entropy conserving 
	formulation at day 12 (right) for the baroclinic instability test case. 
	\blue{Values are given in units of Kelvin, $K$}.}
\label{fig::ce_theta}
\end{center}
\end{figure}

The onset of numerical instability is also observed for the flux form in the vertical component of the
vorticity field at a vertical height of $z\approx 1.5km$ (Fig. \ref{fig::ce_vorticity}), and
the potential temperature at a vertical height of $z\approx 1.25km$ (Fig. \ref{fig::ce_theta}). 
While the flux form solution develops grid scale noise around day 11 (left), this noise 
is not present for the coupled formulation either at day 11 (center) or day 12 (right).

\section{Conclusion}

This article presents energy and entropy conserving spatial formulations
for thermal atmospheric systems in vector invariant form using compatible finite elements. 
Based on the insight that the chain rule is not satisfied for thermodynamic variables on 
discontinuous function spaces, this requirement is negated for the conservation of entropy 
via a novel formulation that bears a close similarity to an averaging of flux form and 
material form transport of thermodynamic variables 
(buoyancy for the thermal shallow water equations and potential temperature for the 3D 
compressible Euler equations), which is also studied. Experiments for the thermal shallow 
water equations show 
that the coupled entropy conserving formulation (an extension of a recent entropy conserving 
discontinuous Galerkin formulation \cite{Ricardo23b}) exhibits superior conservation 
properties to an alternative mixed material/flux form transport formulation, for which 
entropy is not strictly discretely conserved, but for which results are nevertheless stable 
in the absence of any form of stabilisation. Corresponding 
representations of the pressure gradient terms are presented which preserve the skew-symmetric
property, and hence energy conservation, with respect to the thermodynamic transport terms.

Experiments for the 3D compressible Euler equations show that in contrast to a 
non-entropy-conserving potential temperature formulation, the new entropy conserving coupled
formulation is able to run stably without the need for a dissipative term on the thermodynamic
transport and exhibits significantly smaller entropy conservation errors.

\section{CRediT authorship contribution statement}
Kieran Ricardo: Conceptualization, Methodology, Formal analysis, Writing - reviewing and editing.
David Lee: Conceptualization, Formal analysis, Investigation, Methodology, Software, Validation,
Visualization, Writing - original draft.
Kenneth Duru: Methodology, Formal analysis, Writing - reviewing and editing.

\end{document}